\documentclass[11pt,a4paper,reqno]{amsart}
\usepackage[T1]{fontenc}
\usepackage{lmodern}
\usepackage{microtype}

\usepackage[margin=1in]{geometry}
\usepackage{indentfirst}
\usepackage{aliascnt}
\usepackage{amssymb,mathrsfs,needspace,booktabs}
\usepackage{tikz-cd}
\usepackage[all]{xy}
\usepackage{xcolor}
\usepackage[
  colorlinks=true,
  linkcolor=red,
  citecolor=blue,
  urlcolor=blue
]{hyperref}

\newtheorem{maintheorem}{Theorem}

\newtheorem{conj}[maintheorem]{Conjecture}

\newtheorem{theorem}{Theorem}[section]
\newaliascnt{lemma}{theorem}
\newtheorem{lemma}[lemma]{Lemma}
\aliascntresetthe{lemma}

\newaliascnt{proposition}{theorem}
\newtheorem{proposition}[proposition]{Proposition}
\aliascntresetthe{proposition}

\newaliascnt{corollary}{theorem}
\newtheorem{corollary}[corollary]{Corollary}
\aliascntresetthe{corollary}

\theoremstyle{remark}
\newaliascnt{remark}{theorem}
\newtheorem{remark}[remark]{Remark}
\aliascntresetthe{remark}

\numberwithin{equation}{section}

\newcommand{\CC}{\mathbb C}
\newcommand{\PP}{\mathbb P}

\newcommand{\Om}{\Omega}

\DeclareMathOperator{\rk}{rank}

\newcommand{\OO}{\mathcal O}

\renewcommand{\Om}{\Omega_B^1(\log D)}
\newcommand{\KS}{\operatorname{KS}}
\newcommand{\Hom}{\operatorname{Hom}}

\DeclareMathOperator{\gcdop}{gcd}

\allowdisplaybreaks[1]

\hypersetup{pdftitle={Maximal destabilization, Arakelov equality and Kobayashi geodesics},pdfauthor={Tianzhi Hu, Xin Lu, Kang Zuo}}
\title[Maximal destabilization and Arakelov equality]{Maximal destabilization of direct image sheaves,\\ Arakelov equality and Kobayashi geodesics}
\author{Tianzhi Hu}
\address{School of Mathematics and Statistics, Wuhan University, Wuhan, Hubei 430072, P.R. China}
\email{hutianzhi@whu.edu.cn}

\author{Xin L\"u}
\address{School of Mathematical Sciences, Key Laboratory of MEA (Ministry of Education) \& Shanghai Key Laboratory of PMMP, East China Normal University, Shanghai 200241, China}
\email{xlv@math.ecnu.edu.cn}

\author{Kang Zuo}
\address{School of Mathematics and Statistics, Wuhan University, Wuhan, Hubei 430072, P.R. China; Shanghai Institute for Mathematics and Interdisciplinary Sciences (SIMIS), Shanghai, P.R. China; Institut f\"ur Mathematik, Johannes Gutenberg-Universit\"at Mainz, Mainz 55099, Germany}
\email{zuok@uni-mainz.de}
\date{}
\subjclass[2020]{14D06, 14H15, 14J29, 32G15}
\keywords{Pluricanonical sheaf, Arakelov inequality, Teichm\"uller curve, Kobayashi geodesic}

\theoremstyle{plain}
\newaliascnt{conjecture}{theorem}

\aliascntresetthe{conjecture}

\theoremstyle{definition}
\newaliascnt{example}{theorem}

\aliascntresetthe{example}

\DeclareMathOperator{\op}{op}

\begin{document}
\begin{abstract}
Let $f:X\to B$ be a projective semistable family of $n$-folds over a
smooth projective complex curve with connected fibres.
Put $E_\nu=f_*\omega_{X/B}^{\otimes\nu}$ and $d=\deg\omega_B(D)>0$,
where $D$ is the reduced divisor of singular values. We prove the strict Arakelov inequality $\mu(W)<n\nu d/2$ for every subsheaf $W\subseteq E_\nu$ of $\rk W\ge2$. For families of curves of genus at least two, the Arakelov equality for a line bundle $L_\nu\subset E_\nu$ makes $B\setminus D$ a finite orbifold cover of a Teichm\"uller
curve. The line $L_\nu$ is a power of a line in $E_1$ for odd $\nu$ and of a line in $E_2$ for even $\nu$. In higher dimension, the Arakelov equality for a line in $E_1$ shows that the classifying map $\varphi:B\setminus D\to\mathcal M$ is a Kobayashi totally geodesic curve where $\mathcal M$ is the moduli of the fibers of the smooth family $X^\circ:=f^{-1}(B\setminus D)\overset{f}{\to} B\setminus D$. For non-isotrivial fibrations of such curves over $\PP^1$, we deduce that $\kappa(X)\ge1$ requires at least seven singular fibres and determine the graded canonical ring when there are six.  
\end{abstract}
\maketitle

\section{Introduction}
Throughout this paper, we work over the field of complex numbers $\mathbb{C}$.
Let $f:X\to B$ be a projective semistable morphism from a smooth projective $(n+1)$-fold $X$ to a smooth projective curve $B$
with connected fibers of dimension $n\ge1$. 
Let $D\subseteq B$ be the reduced divisor of singular values, and put
\[
 U=B\setminus D,\qquad \Om=\omega_B(D),\qquad
 d=\deg\Om,\qquad E_\nu=f_*\omega_{X/B}^{\otimes\nu}\quad(\nu\ge1).
\]

The positivity of the pluricanonical direct images $E_\nu=f_*\omega_{X/B}^{\otimes\nu}$ plays a central role in the study of Iitaka's conjecture and in the  hyperbolicity of moduli spaces \cite{VZ02}. Viehweg proved the weak positivity of these sheaves \cite[Theorem~III]{Vie83}, and Kawamata established the nefness of $E_\nu$ \cite[Theorem~1]{Kaw82}. In this paper, we study upper bounds for the slopes of subsheaves of $E_\nu$ and the geometry of the equality case.

Throughout, we \textbf{always assume} that $d=\deg\Om>0$. For a nonzero vector bundle $W$ on $B$, write
$\mu(W)=\deg W/\rk W$ and denote its maximal Harder--Narasimhan slope
by $\mu_{\max}(W)$; namely,
$$\mu_{\max}(W):=\big\{\mu(V)~\big|~0\neq V\subseteq W\big\}.$$

It is a very classical problem to bound $\mu(E_{\nu})$ or more generally $\mu(W)$ for any non-zero subbundle $W\subseteq E_{\nu}$
from above by using the information from the base $B$ and the general fiber of $f$ since Parshin-Arakelov's solution to Shafarevich's conjecture over the function field case, cf. \cite{par68,ara71}.
Such kinds of inequalities are nowadays called Arakelov inequalities.
Faltings \cite{fal83} and Deligne \cite{del84} studied the Arakelov inequality for a semistable family
$f:\,X \to B$ of $n$-dimensional abelian varieties;
namely, they showed in this case that
\begin{equation}\label{eq:weak}
\mu(W)\leq \frac{n\nu}2d,\qquad\forall~W\subseteq E_\nu.
\end{equation}
Since then, there have been lots of beautiful works on Arakelov inequalities for families of varieties over curves and also for systems of Hodge bundles over cuvers,
cf. \cite{bcgt,jz02,MVZ06,pet00,VZ04,VZ06}.
It raises therefore a natural question whether the Arakelov equalities can happen.

For a non-isotrivial family of curves of genus $g\geq2$, Tan \cite{Tan95} and Liu \cite{Liu96}
independently proved the strict Arakelov inequality for $W=E_1$,
\begin{equation}\label{eq:LiuTan}
 \mu(E_1)<\frac 12d.
\end{equation}
Viehweg-Zuo \cite{VZ04} gave a complete characterization
of the Arakelov inequalities when
$f:\,X\to B$ is a semistable family of abelian varieties.
Recently, L\"u, Yang and Zuo \cite{LYZ22} proved the strictness for a non-isotrivial family when
$W$ defines a birational morphism over $B$. For families of general
type they conjectured that $\rk W\ge2$ suffices to imply the strictness of the Arakelov inequality in \eqref{eq:weak}; see
\cite[Theorem~1.6 and Conjecture~1.7]{LYZ22}. We prove this rank
criterion for fibres of arbitrary Kodaira dimension.

\begin{maintheorem}[Strictness in higher rank]\label{thm:strict}
For every $\nu\ge1$, a nonzero coherent subsheaf $W\subseteq E_\nu$
satisfies \eqref{eq:weak}. If $\rk W\ge2$, then
\begin{equation}\label{eq:strict}
 \mu(W)<\frac{n\nu}{2}d.
\end{equation}
If $\mu_{\max}(E_\nu)=n\nu d/2$, the maximal destabilizing subbundle
is the unique line bundle of degree $n\nu d/2$ in $E_\nu$.
\end{maintheorem}

Our theorem recovers the strict Arakelov inequality of Tan and Liu
for semistable family of curves with $W=E_\nu$ mentioned above.
%For a non-birationally isotrivial family with general fibre of general
%type, the hypothesis $d>0$ follows from \cite[Theorem~1.1]{LYZ22}.
For the full direct image, Theorem~\ref{thm:strict} also gives the
Kodaira-dimension assertion of \cite[Conjecture~1.4]{LYZ22}. In fact,
it suffices to have equality in one sufficiently large degree.

\begin{corollary}\label{cor:large-pluricanonical-equality}
Let $F$ be the geometric generic fibre of $f$. There is an integer
$\nu_0(F)\ge1$ such that, if for some $\nu\ge\nu_0(F)$ one has
\[
 E_\nu\ne0,\qquad \mu(E_\nu)=\frac{n\nu}{2}d,
\]
then $\kappa(F)=0$.
\end{corollary}

\begin{proof}
Put $P_a(F)=h^0(F,\omega_F^{\otimes a})$ for $a\ge0$.
By Theorem~\ref{thm:strict}, equality forces
$P_\nu(F)=\rk E_\nu=1$.
If $\kappa(F)>0$, choose $m\ge1$ with $P_m(F)\ge2$.
For each residue class $r$ modulo $m$ in which a nonzero
plurigenus occurs, choose a nonnegative degree $a_r$ in that class
with $P_{a_r}(F)>0$. Multiplication by nonzero pluricanonical
forms gives
\[
 P_{a_r+tm}(F)\ge P_m(F)\ge2\qquad(t\ge1).
\]
Thus every nonzero plurigenus in degree
$\nu\ge m+\max_r a_r$ is at least two, which rules out equality
in this range. Since $E_\nu\ne0$ implies $\kappa(F)\ge0$,
the conclusion follows.
\end{proof}

\medskip

In view of Theorem~\ref{thm:strict}, it is natural to ask when the Arakelov equality holds by some line subbundle of $E_\nu$, what special geometric property the family $X/B$ has. For families of abelian varieties, Viehweg and Zuo related Arakelov equality to Shimura curves through the maximal Higgs bundle \cite{VZ04}. For semistable K3 surfaces families with infinite local monodromy, Sun, Tan, and Zuo showed that Arakelov--Yau equality implies that generic fiber has Picard number $19$ and the modular property of the base after a finite cover \cite[Theorem~0.1]{STZ03}.

We now focus on the family of curves of genus $g\ge2$. Damjanovi\'c studied this question and conjectured that if there is a line subbundle $L_\nu\subset E_\nu$ with $\deg L_\nu=\nu d/2$ for some $\nu\ge2$, then the classifying image of $U$ in $\mathcal M_g$ is a \textbf{Teichm\"uller curve}; see \cite[pp.~6--7]{Damjanovic18}.

For the case $\nu=1$, M\"oller proved that the base is a Teichm\"uller curve generated by holomorphic one-forms if and only if the Arakelov equality holds \cite{Moller06}. For $\nu=2$, Hu and Zuo \cite[Theorem~C]{HZ26} proved the base is a Teichm\"uller curve generated by quadratic differential forms if and only if the Arakelov equality holds. Theorem~\ref{thm:characterization} answers Damjanovi\'c's question for every $\nu\ge1$.

\begin{maintheorem}[The Arakelov equality for a family of curves]
\label{thm:characterization}
Assume $n=1$ and the smooth fibres have genus $g\ge2$. For $\nu\geq 1$, if $\mu_{\max}(E_\nu)$ attains its maximum $\nu d/2$, then the maximal destabilizing subbundle is a line bundle $L_\nu\subset E_\nu$ of degree $\nu d/2$ by Theorem~\ref{thm:strict}. Moreover: 
\begin{enumerate}
\item If $\nu$ is odd, there exists a line bundle $L_1\subset E_1$ such that
\[
\deg L_1=\frac d2,\qquad L_\nu=L_1^{\otimes\nu}.
\]

\item If $\nu=2k$ is even, there exists a line bundle $L_2\subset E_2$ such that
\[
\deg L_2=d,\qquad L_\nu=L_2^{\otimes k}.
\]

\item In either case, the image of $U$ in $\mathcal M_g$ is a
Teichm\"uller curve.
\end{enumerate}
\end{maintheorem}

\begin{remark}
Theorem~\ref{thm:characterization} is a \textbf{characterization theorem} of Teichm\"uller curves. Let $C\subset\mathcal M_g$ be a Teichm\"uller curve generated by a holomorphic quadratic differential $q$. After passing to a suitable finite cover of $C$, let $f:X\to B$ be the associated semistable family, with $U=B\setminus D$.

If $q=\omega^2$ for a holomorphic one-form $\omega$, M\"oller proves that there is a line bundle $L_1\subset E_1$ of degree $d/2$ \cite[Lemma~3.1]{Moller06}. Its powers $L_1^{\otimes\nu}\subset E_\nu$ attain the Arakelov bound for every $\nu\ge1$.

If $q$ is not a square, Hu and Zuo prove that $L_2\simeq\Om$ is a line subbundle of $E_2$ \cite[Theorem~C]{HZ26}. Its powers $L_2^{\otimes k}\subset E_{2k}$ attain the Arakelov bound for every $k\ge1$. 

Together with Theorem~\ref{thm:characterization}, the Teichm\"uller curves generated by an Abelian differential can be characterized by the Arakelov equality for all $E_{\nu}$, and the Teichm\"uller curves generated by a quadratic differential can be characterized by the Arakelov equality for all $E_{2k}$.
\end{remark}

\medskip

In higher dimension, the natural generalization of the Teichmuller curves in $\mathcal M_g$ is the Kobayashi totally geodesic curves in the associated moduli space defined in detailed later. M\"oller and Viehweg
characterized Kobayashi totally geodesic curves in the Siegel modular variety $\mathcal A_g$ by the existence of a splitting of the logarithmic Kodaira-Spencer map \cite{MV10}. We aim to prove that the Arakelov equality for a line subbundle of $E_1$ of $X/B$ implies that the associated classifying map has Kobayashi totally geodesic image.

We first recall some terminology. Fix a polynomial $\mathsf h\in\mathbb Q[a,b]$, and let $M_{\mathsf h}$ be Viehweg's quasi-projective coarse moduli scheme of pairs $(Y,\mathcal L)$, where $Y$ is a smooth projective variety with semiample canonical bundle, $\mathcal L$ is an ample line bundle, and $\mathsf h(a,b)=\chi(Y,\mathcal L^{\otimes a}\otimes\omega_Y^{\otimes b})$ \cite[Theorem~1.13]{Vie95}.
For the theorem below, let $\mathcal M$ be a connected fine moduli space for these pairs with a universal polarized family, and let $\varphi:U\to\mathcal M$ be the classifying map of the smooth polarized family $f^{-1}(U)\to U$.

For a connected complex space $S$, the \textbf{Kobayashi pseudodistance} between $x,y\in S$ is defined by
\[
d_S(x,y)=\inf\sum_{i=1}^{m}d_\Delta(a_i,b_i),
\]
where $d_\Delta$ is the Poincar\'e distance on the unit disk $\Delta$ and the infimum is taken over all finite chains of holomorphic maps $\psi_i:\Delta\to S$ with $a_i,b_i\in\Delta$ and $\psi_i(a_i)=x_{i-1}$, $\psi_i(b_i)=x_i$ for points $x=x_0,\ldots,x_m=y$ \cite[Section~1]{MV10}.

By To and Yeung \cite[Theorem~1]{TY15}, a smooth moduli space $\mathcal M$ carrying an effectively parametrized universal family of canonically polarized manifolds is Kobayashi hyperbolic, i.e., $d_{\mathcal M}(x,y)>0$ whenever $x\ne y$.

Following \cite[Definition~1.1]{MV10}, we call $\varphi:U\to\mathcal M$ a \textbf{Kobayashi totally geodesic curve} if its lift $\widetilde\varphi:\Delta\simeq\widetilde U\to\widetilde{\mathcal M}$ satisfies
\[
d_{\widetilde{\mathcal M}}\bigl(\widetilde\varphi(z),\widetilde\varphi(w)\bigr)=d_\Delta(z,w)\qquad\text{for all }z,w\in\Delta.
\]

\begin{maintheorem}[Canonical equality and the Kobayashi geodesic]
\label{thm:geodesic}
Let $f:X\to B$ have relative dimension $n\ge1$, and suppose that $L\subset E_1$ is a line bundle with $\deg L=nd/2$.
Let $\varphi:U\to\mathcal M$ be the  associated classifying map into the moduli space $\mathcal M$. Then the image of $U$ under this map is a Kobayashi totally geodesic curve in $\mathcal M$.
\end{maintheorem}

Theorems~\ref{thm:characterization} and~\ref{thm:geodesic} suggest the following
conjecture. The open case is $n\ge2$ and $\nu\ge2$.

\begin{conj}[Pluricanonical equality and the Kobayashi geodesic]
\label{conj:pluricanonical}
Let $f:X\to B$ have relative dimension $n\ge1$, and suppose that $L_\nu\subset E_\nu$ is a line bundle with $\deg L_\nu=n\nu d/2$.
Let $\varphi:U\to\mathcal M$ be the  associated classifying map into the moduli space $\mathcal M$. Then the image of $U$ under this map is a Kobayashi totally geodesic curve in $\mathcal M$.
\end{conj}

\medskip

We apply the slope inequalities to the number $s$ of singular fibres of a non-isotrivial semistable fibration of curves of genus $g\ge2$ over $\PP^1$. Beauville proved $s\ge4$ \cite{Beauville81}, and Tan improved this bound to $s\ge5$ \cite{Tan95}. If $\kappa(X)\ge0$, Tan, Tu and Zamora proved $s\ge6$ \cite{TTZ05}. 
L\"u, Tan and Zuo proved $s\ge7$ when $X$ is of general type and left open the existence of a fibration with $\kappa(X)=1$ and $s=6$
\cite[Theorems~1 and~3]{LTZ18}.

When $s=6$, the projection formula identifies a nonzero section of $\nu K_X$ with a nonzero map $\OO_{\PP^1}(2\nu)\to E_\nu$. Since
$d=4$, its image attains the Arakelov bound. The uniqueness in
Theorem~\ref{thm:strict} gives $h^0(X,\nu K_X)\le1$ and rules out
$\kappa(X)=1$. The power identities in Theorem~\ref{thm:characterization} then determine the multiplication in the canonical ring.

\begin{maintheorem}\label{thm:intro-surfaces}
For a non-isotrivial semistable curve fibration of fibre genus $g\ge2$ over $\PP^1$, let $s$ be the number of singular fibres.
If $s\le5$, then $\kappa(X)=-\infty$.
If $s=6$, the graded canonical ring of $X$ is isomorphic to one of
\[
 \CC,\qquad \CC[u]\quad(\deg u=1),\qquad
 \CC[v]\quad(\deg v=2).
\]
In particular, $\kappa(X)\ge1$ implies $s\ge7$.
Over an arbitrary base, let $f:X\to B$ be a non-isotrivial semistable fibration of curves of genus $g\ge2$ over a smooth projective curve, with $d>0$.
If $\kappa(X)=1$, then $\mu_{\max}(E_\nu)<\nu d/2$ for every $\nu\ge1$.
\end{maintheorem}

In Section~\ref{sec:strict}, we prove Theorem~\ref{thm:strict} by constructing a pencil on a resolution of a fibre power of $X$ over $B$.
The proof uses Viehweg and Zuo's cyclic-cover construction for pluricanonical sections \cite[Section~4]{VZ02} and the construction of a Higgs subbundle by iterating the Higgs field on the resulting Hodge line \cite[Section~3.1]{LYZ22}. We use both techniques again in subsequent sections. Section~\ref{sec:maximal} uses lifted Beltrami differential formula and the curvature computations of the Hodge metric to prove Theorem~\ref{thm:characterization}. 
Section~\ref{sec:geodesic} proves Theorem~\ref{thm:geodesic} by comparing Hodge distance with Kobayashi distance. The applications to fibred surfaces occupy
Sections~\ref{sec:surfaces}.

\medskip

\noindent\textbf{AI Declaration:} The authors proposed the research questions and some ideas of this paper. ChatGPT raised the pencil construction by sections in Lemma~\ref{lem:two-sections} and proved Proposition~\ref{prop:Arakelovineq}. AI tools were also used for polishing the writing. All
mathematical statements have been independently checked by the authors, who take full
responsibility for the content.

\section{The Arakelov inequality and strictness in rank at least two}
\label{sec:strict}
In this section, we prove Theorem~\ref{thm:strict}.
In Section~\ref{sec:two-sections}, we use the rational map $\eta:X\dashrightarrow\mathbb P_B(W)$ of L\"u, Yang, and Zuo \cite[Theorem~1.6 and Section~3.2]{LYZ22} to construct a pencil of the $2r\nu$-fold fiber product of $X$ over $B$.
In Section~\ref{sec:common}, we recall the Viehweg--Zuo cyclic cover construction and relate the Hodge metric induced by a line bundle attaining the Arakelov equality with the Poincar\'e metric on $U=B\setminus D$.
In Section~\ref{sec:Arakelov}, we prove an Arakelov-type inequality under suitable hypotheses on the Hodge metric, which allows us to deal with the extra singular divisor coming from the pencil cut and the cyclic covering construction.
Finally, in Section~\ref{sec:proofstrict}, we combine these constructions and estimates to prove Theorem~\ref{thm:strict}.

\subsection{Setup and two independent sections}\label{sec:two-sections}

We use the self-fibre product construction of
\cite[(3.5.1)--(3.5.4)]{VZ02}, as also used in
\cite[Section~2.4.2]{LYZ22}.

Let $f:X\to B$ be the semistable family of $n$-folds in the
introduction, and write
\[
 U=B\setminus D,\qquad \Om=\Omega_B^1(\log D),\qquad
 d=\deg\Om>0,\qquad E_\nu=f_*\omega_{X/B}^{\nu}.
\]
By the known Arakelov inequality \eqref{eq:weak}
\cite[Theorem 1.3]{LYZ22}, it remains to exclude the equality for a saturated subsheaf $W\subseteq E_\nu$ of rank
$r=\rk W\ge2$. Suppose, therefore, that $r\geq2$ and
\begin{equation}\label{eq:det-degree}
 \deg W=\frac{rn\nu d}{2}.
\end{equation}

Let $X_B^{2r\nu}:=\underbrace{X\times_B X\times_B\cdots\times_B X}_{2r\nu-\text{times}}$.
Take a smooth projective toroidal resolution
$\pi:X^{[2r\nu]}\to X_B^{2r\nu}$ and let
$f^{[2r\nu]}:X^{[2r\nu]}\to B$ be the induced projection.
Put $\Sigma:=\bigl((f^{[2r\nu]})^{-1}(D)\bigr)_{\mathrm{red}}$,
where $(\cdot)_{\mathrm{red}}$ denotes the associated reduced closed subscheme.
The standard toroidal construction \cite[Chapters I--II]{KKMS}
gives an exact sequence
\begin{equation}\label{eq:logsequence}
 0\longrightarrow (f^{[2r\nu]})^*\Om\longrightarrow\Omega_{X^{[2r\nu]}}^1(\log\Sigma)\longrightarrow\Omega_{X^{[2r\nu]}/B}^1(\log\Sigma)
 \longrightarrow0
\end{equation}
with locally free quotient of rank $2r\nu n$, and
\begin{equation}\label{eq:logcrepant}
 \Lambda:=\det\Omega_{X^{[2r\nu]}/B}^1(\log\Sigma)
 =\omega_{X^{[2r\nu]}/B}\bigl(\Sigma-(f^{[2r\nu]})^*D\bigr)
 =\pi^*\omega_{X_B^{2r\nu}/B}.
\end{equation}
The last equality follows from the logarithmic crepantness of
toroidal modifications; cf.~\cite[Theorems~2.5--2.6]{Ambro16}.
Indeed, if $\mathcal D$ denotes the inverse image of $D\subset B$
on $X_B^{2r\nu}$, then
\[
 \omega_{X^{[2r\nu]}}(\Sigma)
 \simeq \pi^*\bigl(\omega_{X_B^{2r\nu}}(\mathcal D)\bigr).
\]
Semistability gives $\pi^*\mathcal D=(f^{[2r\nu]})^*D$, and the
assertion follows by tensoring with
$(f^{[2r\nu]})^*(\omega_B(D)^{-1})$.

\begin{lemma}\label{lem:two-sections}
There are sections
\begin{equation*}
 t_0,t_1\in H^0\!\left(X^{[2r\nu]},\Lambda^\nu\otimes(f^{[2r\nu]})^*(\det W)^{-2\nu}\right)
\end{equation*}
whose restrictions to the geometric generic fibre of
$f^{[2r\nu]}$ are linearly independent.
\end{lemma}

\begin{proof}
As in \cite[Theorem~1.6 and Section~3.2]{LYZ22}, the evaluation map
\[
\operatorname{ev}:f^*W\longrightarrow\omega_{X/B}^{\nu}
\]
defines a rational map $\eta:X\dashrightarrow\mathbb P_B(W)$.
More precisely, for $x\in X_b$ at which evaluation is nonzero, we have
\[
 \eta_b(x)=[\operatorname{ev}_x],\qquad
 \operatorname{ev}_x:W_b\twoheadrightarrow
 (\omega_{X/B}^{\nu})_x,\qquad s\longmapsto s(x).
\]

For $1\le i\le 2r\nu$, let $q_i:X^{[2r\nu]}\to X$ be the
composition of $\pi$ with the $i$-th projection of $X_B^{2r\nu}$.
For distinct indices $i_1,\ldots,i_r$, let $\Delta_{i_1,\ldots,i_r}$
be the determinant of the map
\[
(f^{[2r\nu]})^*W\xrightarrow{(q_{i_1}^*\operatorname{ev},\cdots, q_{i_r}^*\operatorname{ev})}\bigoplus_{k=1}^r q_{i_k}^*\omega_{X/B}^{\nu}.
\]
Thus $\Delta_{i_1,\ldots,i_r}$ is a global section of $(f^{[2r\nu]})^*(\det W)^{-1}\otimes\bigotimes_{k=1}^r q_{i_k}^*\omega_{X/B}^{\nu}$. In local frames, it is the determinant of homogeneous representatives of $\eta_b(x_{i_1}),\ldots,\eta_b(x_{i_r})$.

Let $\sigma:=(r\;\;r+1)$ be a permutation. Define
\[
 t_0:=\prod_{a=0}^{2\nu-1}\Delta_{ar+1,\ldots,(a+1)r},
 \qquad
 t_1:=\prod_{a=0}^{2\nu-1}\Delta_{\sigma(ar+1),\ldots,\sigma((a+1)r)}.
\]
By \eqref{eq:logcrepant} and the fibre product construction, we have $\bigotimes_{i=1}^{2r\nu}q_i^*\omega_{X/B}^{\nu}=\Lambda^\nu$. Hence $t_0$ and $t_1$ are sections of $\Lambda^\nu\otimes(f^{[2r\nu]})^*(\det W)^{-2\nu}$.

To prove independence, let $b$ be a generic point of $B$ and identify the fibre of $f^{[2r\nu]}$ with $X_b^{2r\nu}$. Since $W\subset f_*\omega^\nu_{X/B}$, there exists $x_1,\cdots,x_r\in X_b$ such that 
\[
W_b\xrightarrow{(\operatorname{ev}_{x_1},\ldots,\operatorname{ev}_{x_r})}\bigoplus_{i=1}^r(\omega_{X/B}^{\nu})_{x_i}
\] 
is an isomorphism. Let $\xi_1\in X_b^{2r\nu}$ repeating 
$(x_1,\ldots,x_r)$ in all $2\nu$ blocks and $\xi_2\in X_b^{2r\nu}$ by exchanging the $r$-th and $(r+1)$-th coordinates of $\xi_1$. Thus
\[t_0(\xi_1)\ne0,\quad t_0(\xi_2)=0,\quad\text{and}\quad t_1(\xi_1)=0,\quad t_1(\xi_2)\ne0.
\]
Thus $t_0$ and $t_1$ are linearly independent.
\end{proof}

\begin{remark}
Suppose $r=2$ and write $P_i=\eta_b(x_i)\in\mathbb P(W_b)$.
After cancelling the common factors, we have
\[
 \frac{t_1}{t_0}
 =\frac{\Delta_{13}\Delta_{24}}{\Delta_{12}\Delta_{34}}.
\]
This is the cross-ratio of four points $(P_1,P_4;P_3,P_2)$.
\end{remark}

\subsection{Cyclic covers and the equality norm}\label{sec:common}

Let $L$ be a holomorphic line bundle on $B$ and let $h$ be a possibly singular Hermitian metric on $L$. On an open set where $h$ is smooth and positive, choose a local coordinate $z$ and a nowhere vanishing holomorphic frame $e$ of $L$. Writing $h(e)=h(e,e)$ for simplicity, define the \textbf{curvature} of $(L,h)$ by
\[
 \Theta_h=-\partial_z\partial_{\bar z}\log h(e).
\]

Put $\Omega=\Omega_B^1(\log D)$ and $U=B\setminus D$.
Since $\deg\Omega>0$, the curve $U$ carries a complete \textbf{Poincar\'e metric} of Gaussian curvature $-4$. Write this metric locally as $\rho(z)^2|dz|^2$ and let $h_\Omega$ be the induced metric on $\Omega|_U=K_U$. Then
\begin{align}\label{eq:Omega}
 h_\Omega(dz)=\rho(z)^{-2},
 \qquad
 \Theta_{h_\Omega} =\partial_z\partial_{\bar z}\log\rho(z)^2=2\rho(z)^2.
\end{align}

Throughout this subsection, let $p_V:V\to B$ be projective, with $V$ smooth and smooth fibres of dimension $N$ over $U$.
Suppose $\Sigma=(p_V^{-1}D)_{\mathrm{red}}$ has simple normal crossings and \eqref{eq:logsequence} has locally free quotient of rank $N$. Put
\[
\Lambda:=\det\Omega_{V/B}^1(\log\Sigma).
\]
Let $L$ be a line bundle on $B$ of degree $Nd/2$. Give $T:=L^2\otimes\Om^{-N}$ a unitary flat metric $h_T$, and set
\begin{equation}\label{eq:ha}
 h_L:=(h_\Omega^N h_T)^{1/2}.
\end{equation}
This defines a metric on $L|_U$.
We first recall the canonical equality case, then treat the pluricanonical equality case by cyclic covers.

\begin{lemma}\label{lem:canonical-equality}
Suppose $p_V$ is semistable and $0\ne s\in H^0(V,\Lambda\otimes p_V^*L^{-1})$.
Then $s$ generates a flat summand of the primitive middle-cohomology VHS of $p_V^{-1}(U)/U$, with graded Higgs bundle
\begin{equation}\label{eq:maxchain}
 H:=\bigoplus_{j=0}^N H_j=\bigoplus_{j=0}^NL\otimes\Om^{-j},\ \text{and Higgs field}\ \tau_j:H_j\xrightarrow{\sim}H_{j+1}\otimes\Om,
\end{equation}
with $0\le j\le N-1.$ Up to one common positive scalar, its Hodge metric is
\begin{equation}\label{eq:chainmetric}
 h_j=\frac{N!\,j!}{(N-j)!}\,h_Lh_\Omega^{-j}.
\end{equation}
On the universal covering disk $\Delta$ of $U$, the Higgs bundle $(\oplus H_j,\oplus \tau_j)$ is the $N$-th symmetric power of the rank-two uniformizing Higgs bundle, up to
a unitary line twist.
\end{lemma}

\begin{proof}
The flat summand and its strictly maximal Higgs field follow from
\cite[Proposition~1.1(c)]{VZ06}.
The symmetric-power description in \cite[Proposition~3.4]{VZ03}, applied after a finite cover admitting a logarithmic theta characteristic if necessary, gives the standard metric \eqref{eq:chainmetric}, up to one common scalar.
\end{proof}

Fix $\nu\ge1$ and a nonzero section
$s\in H^0(V,\Lambda^\nu\otimes p_V^*L^{-\nu})$.
Put $\mathcal M:=\Lambda\otimes p_V^*L^{-1}$ and define the
$\mathcal O_V$-algebra
\[
 \mathcal A_s:=\bigoplus_{i=0}^{\nu-1}\mathcal M^{-i},
\]
whose multiplication is the usual tensor product, followed by
contraction with $s:\mathcal M^{-\nu}\to\mathcal O_V$
whenever the exponents sum to at least $\nu$.
The associated cyclic cover is
\[
 \gamma_s:Y_s:=\operatorname{Spec}_V(\mathcal A_s)
 \longrightarrow V.
\]
It is finite flat of degree $\nu$. Locally, if $\eta$ is a
frame of $\mathcal M$ and $s=\phi\eta^\nu$, then
\[
 Y_s=\operatorname{Spec}
 \bigl(\mathcal O_V[t]/(t^\nu-\phi)\bigr).
\]
The local sections $t\,\gamma_s^*\eta$ glue to a tautological
section $\xi_s$ of $\gamma_s^*\mathcal M$, satisfying
$\xi_s^\nu=\gamma_s^*s$. For $\nu=1$, the cover is the identity.

Normalize $Y_s$, retaining all components, resolve, and
semistably reduce over a finite map $\tau:C\to B$.
Write $p_Z:Z\to C$ for the resulting family and
$\gamma:Z\to V$ for the induced morphism, so that
$p_V\circ\gamma=\tau\circ p_Z$.
Denote the pullback regular root again by $\xi_s$.
Let $D_C$ be a reduced boundary containing all bad values of
$p_Z$ and $\tau^{-1}D$, and put
\[
 \Pi:=(p_Z^{-1}D_C)_{\mathrm{red}},\qquad
 C^\circ:=C\setminus D_C,\qquad
 Z^\circ:=p_Z^{-1}(C^\circ).
\]
Fix a relative polarization and let
\[
 \mathbb V_Z:=
 \bigl(R^N(p_Z|_{Z^\circ})_*\mathbb C\bigr)_{\mathrm{prim}}
\]
be its primitive weight-$N$ VHS.
Let $(\mathcal V_Z,\nabla_Z,F^\bullet)$ be the Deligne
canonical extension to $C$ with its extended Hodge filtration,
and write
\[
 E_Z:=\bigoplus_{j=0}^N E_Z^{N-j,j},\qquad
 E_Z^{N-j,j}:=\operatorname{Gr}_F^{N-j}\mathcal V_Z,
 \qquad \theta_Z:=\operatorname{Gr}_F\nabla_Z.
\]
Thus $\theta_Z$ takes values in
$E_Z\otimes\Omega_C^1(\log D_C)$. The natural map of base
differentials is
\begin{equation}\label{eq:dtau}
 d\tau:\tau^*\Om\to\Omega_C^1(\log D_C).
\end{equation}
Multiplication by $\xi_s$ and pullback of relative top forms give
\[
 p_Z^*\tau^*L\xrightarrow{\ \cdot\xi_s\ }
 \gamma^*\Lambda\longrightarrow
 \Omega_{Z/C}^N(\log\Pi).
\]
By adjunction this induces a nonzero morphism
$u_0:\tau^*L\to E_Z^{N,0}$.

\begin{lemma}
\label{lem:cyclic-maximal}
The Higgs bundle $(E_Z,\theta_Z)$ contains a saturated Higgs
subbundle $H=\bigoplus_{j=0}^N H^{N-j,j}$ with
$H^{N,0}=\operatorname{im}(u_0)$ and for any $j=0,\cdots,N$
\begin{equation}\label{eq:cyclic-maxchain}
H^{N-j,j}\simeq\tau^*(L\otimes\Om^{-j}),\qquad
 \theta_Z|_{H^{N-j,j}}=\mathrm{id}\otimes d\tau.
\end{equation}
It corresponds to an orthogonal flat summand of $\mathbb V_Z$.
\end{lemma}

\begin{proof}
Following \cite[\S~2.2, (2.3)]{LYZ22} and
\cite[\S~4, (4.0.1)]{VZ02}, write $(F,\tau_F)$ for the deformation Higgs bundle of $p_V$, with
\[
F^{p,q}:=R^q(p_V)_*\bigl(\Omega_{V/B}^p(\log\Sigma)\otimes\Lambda^{-1}\bigr)\big/\mathrm{torsion},
 \qquad p+q=N.
\]
Its Higgs fields $\tau_F^{p,q}$ are the connecting morphisms of
the logarithmic exterior-power sequences twisted by $\Lambda^{-1}$.
Set $L'=\tau^*L$ and $\iota=d\tau$.

The Viehweg--Zuo cyclic-cover comparison \cite[Lemma~4.4(i)]{VZ02} and \cite[Proposition~2.4]{LYZ22}, with base change as in \cite[Lemma~4.1]{VZ02} gives maps
\[
 \rho^{p,q}:\tau^*F^{p,q}\longrightarrow
 E_Z^{p,q}\otimes(L')^{-1}
\]
compatible with $\tau_F$ and $\theta_Z$ via $\iota=d\tau$.
The restriction of $\rho^{N,0}\otimes\mathrm{id}_{L'}$
to the unit section of $\tau^*F^{N,0}$ gives $u_0$ defined above by the root form.

Iterating $\tau_F$ on the unit section of $F^{N,0}=(p_V)_*\mathcal O_V$ gives maps $v_j:\Om^{-j}\to F^{N-j,j}$. Pulling them back, tensoring by $L'$, and applying $\rho^{N-j,j}\otimes\mathrm{id}_{L'}$ gives
\begin{equation}\label{eq:iterates}
 u_j:\tau^*(L\otimes\Om^{-j})\longrightarrow E_Z^{N-j,j},\qquad \theta_Z\circ u_j=u_{j+1}\otimes\iota \quad (0\le j<N).
\end{equation}

Let $\ell$ be the last index for which $u_\ell\ne0$. This defines a Higgs subsheaf $H:=\bigoplus_{j=0}^{\ell}\operatorname{im}(u_j)$ of $(E_Z,\theta_Z)$.
Thus the semistability of $(E_Z,\theta_Z)$ and $\deg L=Nd/2$ give
\[
 0\ge\deg H
 =(\deg\tau)\sum_{j=0}^{\ell}(\deg L-jd)
 =\frac{(\deg\tau)(\ell+1)(N-\ell)d}{2}.
\]
Since $d>0$, we have $\ell=N$ and $\deg H=0$.
After taking the saturation of $H$, we get a degree zero Higgs subbundle of $(E_Z,\theta_Z)$ of weight $N$, still denoted by $(H,\theta_Z|_H)$. This proves that $H^{N-j,j}\cong \tau^*(L\otimes\Om^{-j})$.

By the non-abelian Hodge correspondence \cite{Simpson90}, this degree-zero Higgs subbundle corresponds to a flat subbundle. Under the identifications given by the $u_j$, the successive Higgs fields are $\mathrm{id}\otimes d\tau$, as in
\eqref{eq:cyclic-maxchain}. 
\end{proof}

\begin{lemma}\label{lem:common}
For every nonzero section
$s\in H^0(V,\Lambda^\nu\otimes p_V^*L^{-\nu})$, $\nu\ge1$,
there is a constant $c_s>0$ such that
\begin{equation}\label{eq:normrigidity}
I_y(s\ell_0^\nu):=\int_{V_y}|s\ell_0^\nu|^{2/\nu}
 =c_s h_L(\ell_0)(y)
\end{equation}
for any point $y\in U$ and local frame $\ell_0$ of $L$ near $y$.
In particular, $s$ is nonzero on every smooth fibre.
\end{lemma}

The integral uses the positive coordinate volume density:
if $s\ell_0^\nu=\phi(dz_1\wedge\cdots\wedge dz_N)^\nu$,
its density is
$|\phi|^{2/\nu}\prod_j(i/2)dz_j\wedge d\bar z_j$.

\begin{proof}
Let $\bigoplus_{j=0}^N h_{H^{N-j,j}}$ denote the harmonic
Hodge metric on the flat summand $H$ of
Lemma~\ref{lem:cyclic-maximal}.
We prove that, under the identifications
\eqref{eq:cyclic-maxchain},
\[
 h_{H^{N-j,j}}=a\,\tau^*h_j,
 \qquad 0\le j\le N,
\]
for one constant $a>0$, where $h_j$ are the metrics of
Lemma~\ref{lem:canonical-equality}, given by
\eqref{eq:chainmetric}.

Indeed, the metrics $\tau^*h_j$ solve the harmonic equations for the Higgs bundle in \eqref{eq:cyclic-maxchain}.
They are adapted to the full boundary: they extend smoothly at new bad values lying over $U$ and have pullback logarithmic growth over $D$.
Since $H$ is stable, uniqueness of the adapted harmonic metric
\cite[Sections~4 and~6]{Simpson90} therefore proves the claim.
Since $h_0=h_L$, the top Hodge metric is $a\,\tau^*h_L$.

For a point $c\in C^\circ$ and $y=\tau(c)\in B$, the section $u_0(\tau^*\ell_0)$ represents a top form on $Z_c$.
Locally, if $s\ell_0^\nu=\phi(dz_1\wedge\cdots\wedge dz_N)^\nu$, then $t^\nu=\phi$ and
$u_0(\tau^*\ell_0)=t\,\gamma^*(dz_1\wedge\cdots\wedge dz_N)$.
Change of variables, retaining all components of the cover, gives
\[
 \|u_0(\tau^*\ell_0)\|^2
 =\nu\int_{V_y}|s\ell_0^\nu|^{2/\nu}.
\]
Resolution does not change top-form integrals.
The top-metric identity therefore proves
\eqref{eq:normrigidity} in a dense open subset of $U$ after
pulling back to $C$, and hence in that subset of $U$.
Continuity of the fibre integrals extends the identity over
$U$.
\end{proof}

\subsection{The Arakelov-type inequality}\label{sec:Arakelov}

Throughout, a weight-$q$ $\mathbb C$PVHS on a smooth connected
quasi-projective variety $S$ is assumed to be effective: its Hodge bundles satisfy
\[
 E^{q-p,p}=0
 \qquad\text{if }p<0\text{ or }p>q.
\]

A $\mathbb C$PVHS is called \textbf{geometric} if it is a direct summand of
\[
 (R^q f_*\mathbb Q)\otimes_{\mathbb Q}\mathbb C
\]
for some smooth projective morphism $f:Y\to S$.
The summand itself need not be defined over
$\mathbb Q$.

\begin{proposition}\label{prop:Arakelovineq}
Let $U=B\setminus D$ be a punctured Riemann surface and $U_0\subset U$ with $U\setminus U_0$ finite. Let $\mathbb V$ be a geometric $\mathbb C$PVHS of weight $q$ on $U_0$,
and let $\iota:L|_{U_0}\hookrightarrow F^q\mathcal V$ be an inclusion, where $L$ is a line bundle on $B$ and $\mathcal V$ is the holomorphic bundle underlying $\mathbb V$. Suppose that the pullback Hodge metric extends to a \textbf{real analytic metric} $h$ on $L|_U$.
At each point of $D$, suppose that in a local frame $e$ of $L$
and a local parameter $t$, the Hodge norm
\[
 \|\iota(e)\|_{\rm Hod}^2=O\bigl((1+|\log|t||)^k\bigr)
\]
for some $k\geq 0$. Then
\[
 \deg L\leq \frac q2\deg\Omega_B^1(\log D).
\]
\end{proposition}

Let $D_0:=B\setminus U_0$. Under the hypotheses of Proposition~\ref{prop:Arakelovineq}, the standard Arakelov inequality (cf.~\cite[\S~3.1]{LYZ22}), applied to $L|_{U_0}\hookrightarrow F^q\mathcal V$ with boundary $D_0$, gives
\[
 \deg L\leq \frac{q}{2}\deg\Omega_B^1(\log D_0).
\]
Proposition~\ref{prop:Arakelovineq} sharpens this estimate by using the real-analytic extension of the induced Hodge metric across $U\setminus U_0$ and its logarithmic growth near $D$.

\begin{proof}[Proof of Proposition~\ref{prop:Arakelovineq}]
Let $\mathbb J$ on $U_0$ be the smallest flat  subsystem containing the given line.
By \cite[Corollary~4.2.8(iii)(a)]{Deligne}, each
Hodge-torus translate $t\mathbb J$ is flat.
Since $F^q\mathcal V$ is purely of type $(q,0)$,
$t\mathbb J$ still contains the given line.
Minimality and equality of ranks give $t\mathbb J=\mathbb J$.
Thus $\mathbb J$ is a complex subvariation, and its polarization is nondegenerate by Hodge--Riemann positivity.

We first use the real analyticity of $h$ to extend $\mathbb J$
as a polarized variation across $U\setminus U_0$.
Fix $b\in U\setminus U_0$ and choose a disk $\Delta_b$ centered
at $b$ such that $\Delta_b^*=\Delta_b\setminus\{b\}\subset U_0$.
Let $e$ be a holomorphic frame of $L$ over $\Delta_b$, and let
$p_b:\widetilde{\Delta_b^*}\to\Delta_b^*$ be the universal cover.
In a flat trivialization of $p_b^*\mathbb J$, write $\sigma$ for
the pullback of $\iota(e)$. Choose the sign of the flat Hermitian
polarization $Q$ so that it agrees with the Hodge metric on $F^q\mathcal V$. Then
\[
Q\bigl(\sigma(\tilde t),\sigma(\tilde t)\bigr)
=h(e)\bigl(p_b(\tilde t)\bigr)>0\qquad\text{for any}\quad \tilde t\in \widetilde{\Delta_b^*}.
\]

For any $\tilde u\in \widetilde{\Delta_b^*}$, we choose a path from $\tilde u$ to $\tilde t$ in $\widetilde{\Delta_b^*}$ and denote it by $\gamma_{\tilde u,\tilde t}$. By the definition of $\mathbb J$, the parallel transform $(\gamma_{\tilde u,\tilde t})_*(\sigma(\tilde u))\in p_b^*\mathbb J_{\tilde t}$. Hence we may define
\[
 K(\tilde t,\overline{\tilde u}):=Q(\sigma(\tilde t),(\gamma_{\tilde u,\tilde t})_*(\sigma(\tilde u))).
\]
By the real analyticity of $h$ across $b$, the above $K(\tilde t,\overline{\tilde u})$ also descends to a function on $\Delta_b^*\times\Delta_b^*$ and extends as a real analytic function on $\Delta_b\times\Delta_b$, which is still denoted by $K$.

Thus for any $u\in \Delta_b^*$, $K$ is single valued on $\Delta_b^*\times \{u\}$. Let $\gamma_b$ be a loop of $\Delta_b^*$ circled around $b$ and let $M_b:\mathbb J\to \mathbb J$ be its monodromy. By the single-valued property of $K$, it follows that
$(M_b-1)\sigma(t)$ is $Q$-orthogonal to every $\sigma(u)$. The nondegeneracy of $Q$ on $\mathbb J$ gives $M_b|_{\mathbb J}=1$.
By the nilpotent orbit theorem \cite[Theorem~4.9]{Schmid},
$\mathbb J$ extends as a polarized variation across $b$ and
therefore over $U$.

Choose a finite cover $\tau:B'\to B$ that is \'etale over $U$
and makes the local monodromies of $\tau^*\mathbb J$ unipotent
along $D':=(\tau^{-1}D)_{\mathrm{red}}$. Write $U_0'=\tau^{-1}(U_0)$, $L'=\tau^*L$ and $\Omega'=\Omega_{B'}^1(\log D')$. Let $(E,\theta)$ be the logarithmic Higgs bundle associated with the Deligne extension of $\tau^*\mathbb J$. Near a point of $D'$, Schmid's Hodge norm estimates
\cite[Theorem~6.6$'$]{Schmid} show that the Hodge norms of any frame of $E$ have logarithmic growth. 

We extend the pullback of $\iota$ from $U_0'$ to a sheaf inclusion $L'\hookrightarrow E^{q,0}$.
Near a point of $B'\setminus U_0'$, express the image of a local frame of $L'$ in a frame of the Deligne extension adapted to $E^{q,0}$. The coefficient functions are meromorphic, and we aim to prove that they are, in fact, holomorphic. If any of them has a pole in $B'\setminus U_0'$, then its Hodge norm must be unbounded near this pole and is larger than $c\cdot|t|^{-k}$ for some $c,k>0$, which contradicts the assumption of growth on the Hodge norm.

Let $\ell\leq q$ be the largest index for which the iterated Higgs map
\[
 L'\otimes(\Omega')^{-j}\longrightarrow E^{q-j,j}
\]
is nonzero. For $0\leq j\leq\ell$, its image is isomorphic to $L'\otimes(\Omega')^{-j}$. The direct sum of these images is a Higgs subsheaf of $E$. By Simpson's semistability theorem, $(E,\theta)$ is polystable of degree zero, and hence
\[
 0\geq\sum_{j=0}^{\ell}\bigl(\deg L'-j\deg\Omega'\bigr)=(\ell+1)\left(\deg L'-\frac{\ell}{2}\deg\Omega'\right).
\]
As $\tau$ is \'etale over $U$, pullback induces an isomorphism $\tau^*\Omega_B^1(\log D)\simeq\Omega'$. Both $\deg L'$ and $\deg\Omega'$ thus scale by $\deg\tau$. Dividing by this degree
gives
\[
 \deg L\leq\frac{\ell}{2}\deg\Omega_B^1(\log D)\leq\frac{q}{2}\deg\Omega_B^1(\log D).\qedhere
\]
\end{proof}

\subsection{An extremal pencil and the strict inequality}\label{sec:proofstrict}

\begin{proof}[Proof of Theorem~\ref{thm:strict}]
Put $L=(\det W)^{\otimes 2}$ in this proof.
Take the sections $t_0,t_1$ from Lemma~\ref{lem:two-sections} and square them:
\[
 s_i=t_i^2\in H^0\!\left(X^{[2r\nu]},\Lambda^{2\nu}\otimes(f^{[2r\nu]})^*L^{-2\nu}\right),\qquad i=0,1.
\]
They remain independent on the geometric generic fibre: if $(t_1/t_0)^2$ were constant there, $t_1/t_0$ would be constant as well.
The fibres have dimension $N=2r\nu n$, and \eqref{eq:det-degree} gives
\[
 \deg L=r\nu nd.
\]

Let $\psi:=s_1/s_0$ be a rational function on $X^{[2r\nu]}$ and resolve this rational map by a modification $\widetilde\pi:\widetilde X^{[2r\nu]}\to X^{[2r\nu]}$, obtaining a morphism $\widetilde\psi:\widetilde X^{[2r\nu]}\to\mathbb P^1$.
Since
\[
 \Lambda=\omega_{X^{[2r\nu]}/B}\bigl(\Sigma-(f^{[2r\nu]})^*D\bigr)\subseteq\omega_{X^{[2r\nu]}/B},
\]
pullback gives relative $2\nu$-canonical forms on $\widetilde X^{[2r\nu]}$ with twist $(f^{[2r\nu]}\circ\widetilde\pi)^*L^{-2\nu}$.
We use the same notation $s_i$ for these induced forms.

Choose and fix a general $z_0\in\mathbb P^1$.
After removing finitely many base values, we obtain a nonempty Zariski open $U_0\subset U$ and a smooth projective family
\[
 p_{z_0}:Y_{z_0}:=\widetilde\psi^{-1}(z_0)\cap(f^{[2r\nu]}\circ\widetilde\pi)^{-1}(U_0)\longrightarrow U_0,\qquad p_{z_0}:=(f^{[2r\nu]}\circ\widetilde\pi)|_{Y_{z_0}},
\]
of relative dimension $q:=2r\nu n-1$.

The normal bundle of $Y_{z_0}$ in $\widetilde X^{[2r\nu]}$ is the constant line bundle
\[
 N_{Y_{z_0}/\widetilde X^{[2r\nu]}}\simeq \mathcal O_{Y_{z_0}}\otimes T_{z_0}\mathbb P^1.
\]
Choose a normal coordinate $u$ near $z_0$ and let $v_0:=(\partial/\partial u)_{z_0}$.
Restricting $s_0$ to $Y_{z_0}$ and contracting each of its $2\nu$ relative canonical factors with the vector $v_0$, via adjunction, gives a section
\begin{align*}
 \beta\in H^0\!\left(Y_{z_0},\omega_{Y_{z_0}/U_0}^{2\nu}\otimes p_{z_0}^*L^{-2\nu}\right)
\end{align*}
that is nonzero on the general fibre.

For any $b\in U_0$, put $\widetilde X_b^{[2r\nu]}:=(f^{[2r\nu]}\circ\widetilde\pi)^{-1}(b)$ and
\[
 \widetilde\psi_b:=\widetilde\psi|_{\widetilde X_b^{[2r\nu]}}:\widetilde X_b^{[2r\nu]}\longrightarrow\mathbb P^1.
\]
Give $L$ the metric $h_L$ defined by \eqref{eq:ha}.
For a local frame $\ell_0$ of $L$, define a measure on $\mathbb P^1$ by
\[
 \mu_b(A):=\frac{1}{h_L(\ell_0)(b)}\int_{\widetilde\psi_b^{-1}(A)}|s_0\ell_0^{2\nu}|^{1/\nu},
\]
for any Borel subset $A\subset\mathbb P^1$.

Top-form integrals are unchanged by the modification $\widetilde\pi$.
Since $s_0$ vanishes on $\widetilde\psi_b^{-1}(\infty)$, we have $\mu_b(\{\infty\})=0$.
For every $(a_0,a_1)\in\mathbb C^2$,
\[
 \int_{\mathbb C}|a_0+a_1z|^{1/\nu}\,d\mu_b(z)=\frac{I_b\bigl((a_0s_0+a_1s_1)\ell_0^{2\nu}\bigr)}{h_L(\ell_0)(b)}.
\]
By Lemma~\ref{lem:common}, these integrals are finite and independent of $b\in U_0$.
Since $1/\nu\in(0,2)$, Rudin's equimeasurability theorem \cite[Theorem~1]{Rudin} implies that the measures $\mu_b$ define the same measure $\mu$ on $\mathbb C$ and hence on $\mathbb P^1$.

For a coordinate $(v_1,\cdots,v_q,u)$ near $Y_{z_0}$ with normal coordinate $u$, write
\[
 u=\widetilde\psi,\qquad s_0\ell_0^{2\nu}=\phi(u,v,b)(du\wedge dv_1\wedge\cdots\wedge dv_q)^{2\nu}.
\]
The adjunction construction gives
\[
 \beta\ell_0^{2\nu}=\phi(z_0,v,b)(dv_1\wedge\cdots\wedge dv_q)^{2\nu}.
\]
By Fubini's theorem, the density of $\mu_b$ at $z_0$ is
\[
 \frac{1}{h_L(\ell_0)(b)}\int_{(Y_{z_0})_b}|\beta\ell_0^{2\nu}|^{1/\nu}.
\]
These densities are continuous at regular values.
The measure equality $\mu_b=\mu$ therefore implies, after shrinking $U_0$ if necessary,
\begin{align}\label{eq:metriccompatible}
\int_{(Y_{z_0})_b}|\beta\ell_0^{2\nu}|^{1/\nu}=c_{z_0}h_L(\ell_0)(b),\qquad c_{z_0}>0.
\end{align}
Positivity follows from the nonvanishing of $\beta$ on a reference smooth fibre.

Take the $2\nu$-th root cover of $\beta$ in $\omega_{Y_{z_0}/U_0}\otimes p_{z_0}^*L^{-1}$, normalize it, and resolve it.
After shrinking $U_0$ again, write $g_\beta:Z_\beta\to U_0$ for the resulting smooth projective family and $\mathbb V_\beta$ for its weight-$q$ PVHS for the middle cohomology.
The root form gives a map $\iota:L|_{U_0}\to F^q\mathcal V_\beta$.
By change of variables,
\[
 \|\iota(\ell_0)(b)\|_{\mathrm{Hod}}^2=2\nu\int_{(Y_{z_0})_b}|\beta\ell_0^{2\nu}|^{1/\nu}=2\nu c_{z_0}\,h_L(\ell_0)(b).
\]
Thus $\iota$ is a line inclusion, and its induced metric extends to the positive real analytic metric $2\nu c_{z_0}h_L$ on $L|_U$ by \eqref{eq:metriccompatible}.
Near $D$, the hyperbolic metric satisfies $h_\Omega(dt/t)=O((1+|\log|t||)^2)$.
The unitary factor in \eqref{eq:ha} is smooth across $D$, so the squared Hodge norm of a local frame of $L$ is $O((1+|\log|t||)^N)$.
Proposition~\ref{prop:Arakelovineq} therefore gives
\[
 \frac N2d=\deg L\leq\frac q2d=\frac{N-1}{2}d,
\]
contradicting $d>0$.
This proves strictness in rank at least two.

Finally, if $\mu_{\max}(E_\nu)=n\nu d/2$, the above strictness forces the first Harder--Narasimhan term to have rank one.
Let $W$ be a line subbundle of $E_\nu$ with degree $n\nu d/2$.
For any other line subbundle $W'$ with the same degree, consider the map
\[
 W'\hookrightarrow E_\nu\to E_\nu/W.
\]
If $W'\ne W$, this map is nonzero.
But $\deg W'\leq\mu_{\max}(E_\nu/W)<\mu_{\max}(E_\nu)$, which contradicts the degree assumption on $W'$.
\end{proof}

\begin{corollary}\label{coro:Hom}
If $\deg W=n\nu d/2$ for a line bundle $W$ on $B$, then $\dim\Hom(W,E_\nu)\le1$.
\end{corollary}

\section{The equality case for families of curves}\label{sec:maximal}
For family of curves, we study the Arakelov equality and prove Theorem~\ref{thm:characterization} in this section. In Section~\ref{sec:area}, we prove that for the cyclic covering of $Z\to X$ of the inclusion $s:L_\nu\to E_\nu$, the Beltrami differential $\mu_X$ of a smooth fiber of $X/B$ pulls back to the Beltrami differential of $Z$ in Lemma~\ref{lem:root-deformation}. Then we use the cyclic covering comparison again by Viehweg-Zuo to derive the comparison between the Griffiths curvature of the root bundle of $L_\nu$ and the pullback Teichm\"uller metric of the base $U$. This also generalizes to any subbundle $W\subset E_\nu$ in Proposition~\ref{prop:metric}. When $L_\nu\subset E_\nu$ attains the Arakelov equality, this maximality and the above curvature formula shows that the pullback Teichm\"uller metric of $U$ is exactly the Poincar\'e metric of $U$, which completes the proof.

The curvature estimates Proposition~\ref{prop:metric} used in the proof of Theorem~\ref{thm:characterization} also yield a refined Arakelov inequality for families of curves of genus $g\geq2$ (Theorem~\ref{thm:energy}).
\subsection{Beltrami differential of the cyclic covering and the curvature formulas}
\label{sec:area}

Under the hypotheses of Theorem~\ref{thm:characterization}, we have $n=1$, $g\ge2$, and $\mu_{\max}(E_\nu)=\nu d/2$ for some $\nu\ge1$.
By Theorem~\ref{thm:strict}, the maximal destabilizing subbundle is a line $L_\nu\subset E_\nu$ of degree $\nu d/2$.

Fix the inclusion $s:L_\nu\hookrightarrow E_\nu$.
If $d$ is even, set $D_\nu=0$.
If $d$ is odd, then $\nu$ is even and $D\ne\varnothing$, so choose $b_*\in D$ and set $D_\nu=(\nu/2)b_*$.
In both cases $\deg(L_\nu(-D_\nu))$ is divisible by $\nu$.
Therefore there is a line bundle $A$ on $B$ with $A^{\otimes\nu}\simeq L_\nu(-D_\nu)$.
Let $\jmath:A^{\otimes\nu}\hookrightarrow L_\nu$ be the resulting inclusion, which is an isomorphism over $U$.
By adjunction, $s\circ\jmath$ gives a regular section
\[
 \sigma\in H^0(X,\mathcal M^{\otimes\nu}),\qquad \mathcal M:=\omega_{X/B}\otimes f^*A^{-1}.
\]
As in the cyclic-cover construction of \cite[Section~4]{VZ02} and \cite[Sections~2.4.2 and~3.1]{LYZ22}, let $\pi:Z\to X$ be the normalized cyclic cover associated with $(\mathcal M,\sigma)$, retaining all components, and let $\xi$ be its tautological root, with $\xi^{\otimes\nu}=\pi^*\sigma$.

Put $g:=f\circ\pi$.
After deleting finitely many points of $U=B\setminus D$, we obtain a nonempty Zariski open set $U^\circ\subset U$ on which the zero multiplicities are locally constant. Then $g^\circ:Z^\circ:=g^{-1}(U^\circ)\to U^\circ$ is smooth and projective.
pullback of relative differentials sends the tautological root $\xi$ to a section of $\omega_{Z^\circ/U^\circ}\otimes(g^\circ)^*A^{-1}$, giving a morphism
\begin{align}\label{eq:u0}u_0:A|_{U^\circ}\to g^\circ_*\omega_{Z^\circ/U^\circ}.\end{align}

For a local frame $e$ of $A$ on $U^\circ$, write $s_t=s(\jmath(e^{\otimes\nu}))|_{X_t}$, $\alpha_t=u_0(e)|_{Z_t}$ and $\pi_t=\pi|_{Z_t}$.
Then
\[
 \alpha_t\in H^0(Z_t,K_{Z_t}),\qquad \alpha_t^{\otimes\nu}=\pi_t^*s_t.
\]

Fix a point $0\in U^\circ$ and let $t$ be a local coordinate near this point.
Write $\pi_0:Z_0\to X_0$ for the corresponding fibre map and $\alpha_0$ for its tautological one-form.
Choose smooth Beltrami representatives
\[
 \mu_X\in\mathcal A^{0,1}(X_0,T_{X_0})\qquad \text{and}\qquad \mu_Z\in \mathcal A^{0,1}(Z_0,T_{Z_0})
\]
of the Kodaira--Spencer classes of the two families in the direction $\partial_t$.
Let $\mathrm{Ram}(\pi_0)\subset Z_0$ be the set of all ramification points of $\pi_0$.
For $\mu_X$, define its lift on $Z_0\setminus\mathrm{Ram}(\pi_0)$ by
\[
 \widetilde\mu_X(w) =\mu_X(\pi_0(w))\frac{\overline{\pi_0'(w)}}{\pi_0'(w)}\,d\bar w\otimes\frac{\partial}{\partial w}\in \mathcal A^{0,1}(Z_0\setminus\mathrm{Ram}(\pi_0),T_{Z_0}).
\]
Contraction with $\alpha_0$ extends $\widetilde\mu_X\lrcorner\alpha_0$ to a smooth $(0,1)$-form on $Z_0$.
Let $P$ be harmonic projection onto $H^{0,1}(Z_0)$.
Let $\pi_0^*\mu_X$ be the $\pi_0^*T_{X_0}$-valued $(0,1)$-form on $Z_0$ obtained by pulling back $\mu_X$, with local expression
\[
 \pi_0^*\mu_X =\mu_X(\pi_0(w))\,\overline{\pi_0'(w)}\,d\bar w \otimes \left.\frac{\partial}{\partial z}\right|_{\pi_0(w)}\in\mathcal A^{0,1}(Z_0,\pi_0^*T_{X_0}).
\]

\begin{lemma}
\label{lem:root-deformation}
For smooth representatives $\mu_X$ and $\mu_Z$, there exists $V\in C^\infty(Z_0,\pi_0^*T_{X_0})$ such that
\begin{equation}\label{eq:cover-deformation}
 d\pi_0\,\mu_Z-\pi_0^*\mu_X=\bar\partial V \qquad\text{in }\ \mathcal A^{0,1}(Z_0,\pi_0^*T_{X_0}).
\end{equation}
Moreover, we have
\begin{equation}\label{eq:root-projection}
 P(\mu_Z\lrcorner\alpha_0)=P(\widetilde\mu_X\lrcorner \alpha_0).
\end{equation}
\end{lemma}

\begin{proof}
By the Dolbeault description of the Kodaira--Spencer class \cite[Chapter~I, \S4]{Manetti04}, we may choose smooth lifts $A_X\in C^\infty(X_0,T_X|_{X_0})$ and $A_Z\in C^\infty(Z_0,T_Z|_{Z_0})$ of $\partial_t|_0$ with $\bar\partial A_X=\mu_X$ and $\bar\partial A_Z=\mu_Z$.
Since $\pi$ is a map over the base, the sections $d\pi(A_Z)$ and $\pi_0^*A_X$ project to the same base vector $\partial_t|_0$.
Their difference
\[
 V=d\pi(A_Z)-\pi_0^*A_X
\]
is therefore a smooth section of $\pi_0^*T_{X_0}$.
As $d\pi$ is holomorphic and restricts to $d\pi_0$ on vertical tangent vectors along $Z_0$, we obtain
\[
 \bar\partial V =d\pi_0\,\mu_Z-\pi_0^*\mu_X,
\]
which proves \eqref{eq:cover-deformation}.

Let $x\in X_0$ be a zero of $s_0$ of order $m$.
In a local coordinate $z$ centered at $x$, write $s_0=z^ma(z)(dz)^\nu$, where $a(0)\ne0$, and set
\[
 \nu_0:=\frac{\nu}{\gcdop(\nu,m)},\qquad m_0:=\frac{m}{\gcdop(\nu,m)}.
\]
For each $y\in\pi_0^{-1}(x)$, there is a local coordinate $w$ on $Z_0$ centered at $y$ in which
\begin{equation}\label{eq:rootmodel}
 z\circ\pi_0=w^{\nu_0},\qquad \alpha_0=\nu_0\,w^{\nu_0+m_0-1}u(w^{\nu_0})\,dw, \qquad u^\nu=a.
\end{equation}
We view $d\pi_0$ as a holomorphic section of $K_{Z_0}\otimes(\pi_0^*K_{X_0})^{-1}$.
Thus by \eqref{eq:rootmodel}, $\alpha_0/d\pi_0$ extends holomorphically across $\mathrm{Ram}(\pi_0)$ as a section of $\pi_0^*K_{X_0}$, with local coefficient $w^{m_0}u(w^{\nu_0})$.

On $Z_0\setminus\mathrm{Ram}(\pi_0)$, \eqref{eq:cover-deformation} gives
\[
 \mu_Z\lrcorner\alpha_0-\widetilde\mu_X\lrcorner\alpha_0 =\bar\partial\bigl(V\lrcorner(\alpha_0/d\pi_0)\bigr).
\]
For smooth representatives, both sides extend smoothly to $Z_0$, so the equality holds globally.

For any $\eta\in H^0(Z_0,K_{Z_0})$, the contraction $V\lrcorner(\alpha_0\eta/d\pi_0)$ is therefore a smooth one-form on $Z_0$.
Stokes's theorem yields
\[
 \int_{Z_0}\bigr((\mu_Z-\widetilde\mu_X)\lrcorner\alpha_0\bigr)\eta =\int_{Z_0}\bar\partial \left(V\lrcorner(\frac{\alpha_0\eta}{d\pi_0})\right)=0.
\]
This proves \eqref{eq:root-projection}.
\end{proof}

For a local generator $s_t$ of $L_\nu$, define $I(s_t)=\int_{X_t}|s_t|^{2/\nu}$ and give $L_\nu$ the Hermitian metric $h(s_t)=I(s_t)^\nu$.
With the Hodge norm $\|\alpha_t\|^2=\frac i2\int_{Z_t}\alpha_t\wedge\overline{\alpha_t}$, the cyclic-cover norm calculation in the proof of Lemma~\ref{lem:common} gives
\begin{equation}\label{eq:root-norm}
 \|\alpha_t\|^2=\nu I(s_t),\qquad h(s_t)=\left(\frac{\|\alpha_t\|^2}{\nu}\right)^\nu.
\end{equation}

The weight-one variation $R^1g^\circ_*\CC$ on $U^\circ$ has Hodge bundles
\[
 \mathcal H^{1,0}=g^\circ_*\omega_{Z^\circ/U^\circ},\qquad \mathcal H^{0,1}=R^1g^\circ_*\OO_{Z^\circ},
\]
and Higgs field $\theta:\mathcal H^{1,0}\to\mathcal H^{0,1}\otimes\Omega_{U^\circ}^1$.
Let $D_{\rm Hod}$ be the Chern connection of the Hodge metric on $\mathcal H^{1,0}$.
The associated smooth connection of the Gauss--Manin connection satisfies
\[
 \nabla_{\partial_t}u_0(e) =D_{{\rm Hod},\partial_t}u_0(e)+\theta_{\partial_t}u_0(e).
\]
The Higgs field is the harmonic Kodaira--Spencer contraction, so Lemma~\ref{lem:root-deformation} yields
\begin{align}\label{eq:thetaalpha}
 \theta_{\partial_t}\alpha_t =P(\mu_Z\alpha_t)=P(\widetilde\mu_X\alpha_t).
\end{align}
We define the curvature of $(L_\nu,h)$ by $\Theta_h(t):=-\partial_t\partial_{\bar t}\log h(s_t).$
The Hodge curvature formula and the subbundle curvature formula give the following curvature formula of $(L_\nu,h)$.

\begin{lemma}\label{lem:curvature}
Let $\mathcal A:=\operatorname{im}(u_0)\subset\mathcal H^{1,0}$ be the line defined by \eqref{eq:u0}.
With $p_{\mathcal A}$ denoting orthogonal projection onto $\mathcal A$, define its second fundamental form in $\mathcal H^{1,0}$ by
\[
 \sigma_{\mathcal A} =(1-p_{\mathcal A})D_{\rm Hod}^{1,0}|_{\mathcal A}.
\]
For every Beltrami representative $\mu_X$ of the deformation of $X_t$, one has
\begin{equation}\label{eq:Hodgecurv}
 \Theta_h =\frac{\nu}{\|\alpha_t\|^2} \left(\|\theta_{\partial_t}\alpha_t\|^2 -\|\sigma_{\mathcal A,\partial_t}\alpha_t\|^2\right).
\end{equation}
In particular,
\begin{equation}\label{eq:curvbound}
 \Theta_h \le\nu\frac{\|P(\widetilde\mu_X\alpha_t)\|^2}{\|\alpha_t\|^2} \le\nu\frac{\|\widetilde\mu_X\alpha_t\|^2}{\|\alpha_t\|^2} \le\nu\|\mu_X\|_\infty^2.
\end{equation}
\end{lemma}

\begin{proof}
The curvature of the top Hodge bundle $\mathcal H^{1,0}$ in the direction $\partial_t$ and $\alpha_t$ is $\|\theta_{\partial_t}\alpha_t\|^2$.
The subbundle curvature formula therefore gives
\[
 -\partial_t\partial_{\bar t}\log\|\alpha_t\|^2 =\frac{\|\theta_{\partial_t}\alpha_t\|^2 -\|\sigma_{\mathcal A,\partial_t}\alpha_t\|^2} {\|\alpha_t\|^2}.
\]
Using \eqref{eq:root-norm}, we get \eqref{eq:Hodgecurv}.

Finally, orthogonal projection decreases the $L^2$ norm, and $|\widetilde\mu_X|\le\|\mu_X\|_\infty$ almost everywhere.
Hence by \eqref{eq:thetaalpha}, we have
\[
 \|\theta_{\partial_t}\alpha_t\|=\|P(\widetilde\mu_X\alpha_t)\| \le\|\widetilde\mu_X\alpha_t\| \le\|\mu_X\|_\infty\|\alpha_t\|,
\]
which proves \eqref{eq:curvbound}.
\end{proof}

Let $h_\Omega$ be the Poincar\'e on $U$ with $h_\Omega(dt)=\rho(t)^{-2}$. By \cite{Royden70}, the Teichm\"uller metric of the Kodaira--Spencer class is the infimum of the $L^\infty$ norms of its Beltrami representatives.
Thus the pullback Teichm\"uller metric on $U$ is
\begin{align}\label{eq:pullbackTeichmuller}
\rho_T:=\|\KS_f(\partial_t)\|_T=\inf_{\mu_X}\|\mu_X\|_\infty.
\end{align}
In particular, Lemma~\ref{lem:curvature} gives $\Theta_h\le\nu \rho_T^2$.
Since this estimate does not require slope equality, it also applies to any line subbundle of $E_\nu$.

\bigskip

We now use the local curvature estimate to prove global positivity on the projectivization of a subbundle of $E_\nu$.
Let $w:W\hookrightarrow E_\nu$ be the inclusion of a subbundle, and Let $p:\PP(W)\to B$ be the projective bundle parametrizing all lines of $W$.
Put
\[
 M:=\OO_{\PP(W)}(2)\otimes p^*\Om^{\otimes\nu}.
\]
Let $\sigma$ be a local frame of $\OO_{\PP(W)}(-1)$ over $U_0$.
At a point $\xi_t$ of $\PP(W)$ over $t\in U_0$, let $w_\sigma(\xi_t)$ be a local section of $E_\nu$ defined by $w_\sigma(\xi_t):=w_t(\sigma)\in H^0(X_t,K_{X_t}^{\otimes\nu})$.
Define the area metric $h_W$ on $\OO_{\PP(W)}(-1)$ by
\[
I(w_\sigma):=\int_{X_t}|w_\sigma(t)|^{2/\nu},\qquad h_W(\sigma):=I(w_\sigma)^\nu.
\]
The area metric $h_W$ and the Poincar\'e metric of $\Omega$ define a twisted metric $h_M:=h_\Omega^{\nu}\cdot h_W^{-2}$ on $M$. For a frame $\sigma^{-2}\otimes p^*(dt)^{\otimes\nu}$, its local potential function and curvature current are
\begin{equation}\label{eq:weight}
 u=2\log h_W(\sigma)-\nu p^*\log h_\Omega(dt),\qquad \Theta_M=\frac{i}{2\pi}\partial\bar\partial u.
\end{equation}

\begin{proposition}
\label{prop:metric}
The metric $h_M$ extends to a singular Hermitian metric on $M$ with positive curvature current $\Theta_M$.
Its local potential function is continuous on $\PP(W)|_U$ and smooth on a dense Zariski open subset $G\subset\PP(W)|_U$.
Set
\[
 \gamma:=\frac{i}{\pi}(\rho^2-\rho_T^2)\,dt\wedge d\bar t.
\]
Then
\begin{equation}\label{eq:energy-current}
 \Theta_M\ge\nu p^*\gamma\qquad\text{on }G.
\end{equation}
For every real smooth $(1,1)$-form $0\le\eta\le\gamma$ compactly supported in $U$, the class $c_1(M)-\nu p^*[\eta]$ is nef.
In particular, $M$ is nef.
\end{proposition}

\begin{proof}
By the definition of $w_\sigma$, there is a dense Zariski open subset $G$, on which the zeros of $w_\sigma$ are locally given by pairwise disjoint holomorphic sections, each with constant multiplicity. 
For a holomorphic curve $C\subset G$ transverse to the fibres of $p$, let $\tau$ be the local coordinate with $p_*(\partial_\tau)=\partial_t$. Apply Lemma~\ref{lem:curvature} to the cyclic covering family induced by $w_\sigma$ along $C$.
Its Kodaira--Spencer class is $\KS_f(\partial_t)$, so the curvature $\Theta_{h_W}=-\partial_\tau\partial_{\bar\tau}\log h_W(\sigma)$ satisfies $\Theta_{h_W}\le\nu\rho_T^2$.
By \eqref{eq:Omega} and \eqref{eq:weight},
\[
 \partial_\tau\partial_{\bar\tau}u\ge 2\nu(\rho^2-\rho_T^2).
\]
This gives \eqref{eq:energy-current}.

Since $w_\sigma(\xi_t)$ is nonzero, $I(w_\sigma)$ is well-defined and positive on $U=B\setminus D$ and $u$ extends continuously and plurisubharmonically across $\PP(W)|_U\setminus G$.
By \eqref{eq:energy-current}, the same extension applies to $u$ minus a local potential of $\nu p^*\eta$.

Near a singular fibre, let $t$ be a local parameter on $B$.
In a plumbing chart $xy=t$ of $X/B$ near the boundary, the form induced by $\sigma$ is $A(x,y,t)(dx/x)^\nu$, with $A$ locally bounded.
Each plumbing annulus contributes at most $O(\log(R/|t|))$ for some $R>0$ to $I(w_\sigma)$, and the contribution away from the nodes is locally bounded.
Thus, for some $R>0$,
\begin{equation}\label{eq:boundary}
 I(w_\sigma)= O(\log(R/|t|)),\qquad h_\Omega(dt/t)\ge4\log^2(R/|t|).
\end{equation}
In the frame $\sigma^{-2}\otimes p^*(dt/t)^{\otimes\nu}$ of $M$, the potential function is $u=2\log h_W(\sigma)-\nu p^*\log h_\Omega(dt/t)$.
Since $h_W(\sigma)=I(w_\sigma)^\nu$, the bounds in \eqref{eq:boundary} make this $u$ locally bounded above, so it extends plurisubharmonically across $\PP(W)|_D$.
Let $\varphi_\eta$ be a local potential for $\eta$, i.e. $\frac{i}{2\pi}\partial\bar\partial\varphi_\eta=\eta$.
Since $\eta$ vanishes near $D$, we may take $\varphi_\eta=0$ there, so $u-\nu p^*\varphi_\eta$ also extends plurisubharmonically across $\PP(W)|_D$. The extended potential functions give positive currents $\Theta_M$ and $\Theta_M-\nu p^*\eta$ on $\PP(W)$, representing $c_1(M)$ and $c_1(M)-\nu p^*[\eta]$, respectively.

The curve criterion proves that $M$ is nef.
\end{proof}

\subsection{Proof of Theorem~\ref{thm:characterization}}
We first give the following criterion lemma on when a line subbundle of $E_{k\nu}$ is a $k$-th power of some other line subbundle of $E_\nu$.
\begin{lemma}\label{lem:saturation}
\begin{enumerate}
\item If $R\subset E_\nu$ is a line subbundle, then the power map $R^{\otimes k}\to E_{k\nu}$ is an inclusion for every $k\ge1$.
\item Conversely, let $L\subset E_{k\nu}$ be a line subbundle. If for general $t\in B$, any element $\ell_t\in L_t$ is a $k$-th power of a section in $H^0(X_t,K_{X_t}^\nu)$, then $L$ is the $k$-th power of a unique line subbundle $R\subset E_\nu$.
\end{enumerate}
\end{lemma}

\begin{proof}
For (1), since $R$ has rank one, one can verify directly that any $k$-th power of a non-zero element in $R$ is also non-zero.

For (2), let $\eta=\operatorname{Spec}F$ be the generic point of $B$, where $F=\mathbb C(B)$.
Choose an algebraic closure $\bar F$ of $F$ and let $\bar\eta=\operatorname{Spec}\bar F\to\eta\to B$ be the associated geometric generic point.
Set $X_\eta=X\times_B\eta$, $X_{\bar\eta}=X_\eta\times_\eta\bar\eta$, and $L_{\bar\eta}=L_\eta\otimes_F\bar F$.
By flat base change, we have 
\[
 (E_j)_\eta\cong H^0(X_\eta,\omega_{X_\eta}^{\otimes j}),\quad\text{and}\quad  (E_j)_\eta\otimes_F\bar F \cong H^0(X_{\bar\eta},\omega_{X_{\bar\eta}}^{\otimes j}),\qquad j=\nu,k\nu.
\]

Over the smooth locus $U\subset B$, taking the $k$-th powers defines a proper morphism
\[
 \Phi:\mathbb P(E_\nu|_U)\longrightarrow\mathbb P(E_{k\nu}|_U),\qquad [q]\longmapsto[q^k],
\]
where the projective bundles parametrize lines.
By the assumption, for general point $t\in U$, the point $[L_t]\in \mathbb P(E_{k\nu}|_U)$ lies in the image of $\Phi$, so its generic value $[L_\eta]$ also lies in that image.
After base changing to $\bar F$, there is a $\bar F$-point $q\in H^0(X_{\bar\eta},\omega_{X_{\bar\eta}}^{\otimes\nu})$ satisfying $\Phi_{\bar\eta}([q])=[L_{\bar\eta}]$.

The root line $\bar Fq$ is unique. Indeed, if $q_1^k$ and $q_2^k$ span $L_{\bar\eta}$, then $(q_1/q_2)^k\in\bar F^\times$.
Since $X_{\bar\eta}$ is integral and $\bar F$ is algebraically closed, the ratio $q_1/q_2$ belongs to $\bar F^\times$.
Both the power map and $L_\eta$ are defined over $F$, so every element of $\operatorname{Gal}(\bar F/F)$ preserves $\bar Fq$.
In an $F$-basis $\{e_i\}_{i=1}^m$ of $(E_\nu)_\eta$, we may write $q=\sum a_ie_i$ and assume that one of $a_i$ is $1$ after multiplying $q$ by a non-zero element in $\bar F$.
Every element of $\operatorname{Gal}(\bar F/F)$ then fixes $q$, since it preserves $\bar Fq$.
Thus $R_\eta:=Fq$ is a line in $(E_\nu)_\eta$ whose $k$-th power spans $L_\eta$.

Let $R\subset E_\nu$ be the saturated extension of $R_\eta$.
By (1), the image of $R^{\otimes k}\to E_{k\nu}$ is saturated.
It has generic fibre $L_\eta$, so it equals $L$.
Any other line subbundle with the same property has generic fibre $R_\eta$ and hence equals $R$, which proves the uniqueness.
\end{proof}

\begin{proof}[Proof of Theorem~\ref{thm:characterization}]
Apply Proposition~\ref{prop:metric} to $W=L_\nu$.
The positive curvature current on $\Om^{\otimes\nu}\otimes L_\nu^{-2}$ vanishes because this line bundle has degree zero.
Thus, on the subset $U^\circ$ defined in Section~\ref{sec:area}, \eqref{eq:weight} gives
\begin{align*}
 \Theta_h=\nu\rho^2.
\end{align*}
Choose a Beltrami representative $\mu_X$ of minimum $L^\infty$ norm for $\KS_f(\partial_t)$.
By \eqref{eq:curvbound},
\begin{equation}\label{eq:extremal-equality-chain}
 \rho^2(t) \le\frac{\|P(\widetilde\mu_X\lrcorner\alpha_t)\|^2}{\|\alpha_t\|^2} \le\frac{\|\widetilde\mu_X\lrcorner\alpha_t\|^2}{\|\alpha_t\|^2} \le\rho_T^2(t)\le\rho^2(t),
\end{equation}
where the last inequality follows from the Schwarz lemma of the Poincar\'e metric. Hence both inequalities in \eqref{eq:extremal-equality-chain} are equalities.

Equality in \eqref{eq:Hodgecurv} gives $\sigma_{\mathcal A,\partial_t}\alpha_t=0$.
The equality $\|P(\widetilde\mu_X\lrcorner\alpha_t)\|=\|\widetilde\mu_X\lrcorner\alpha_t\|$ gives $\widetilde\mu_X\lrcorner\alpha_t=\overline{\eta_t}$ with $\eta_t$ holomorphic.
The equality $\|\widetilde\mu_X\lrcorner\alpha_t\|=\|\alpha_t\|\rho(t)$ gives $|\eta_t|=\rho(t)|\alpha_t|$ almost everywhere.
The meromorphic ratio $\eta_t/\alpha_t$ has constant nonzero modulus, so it is constant.
Hence
\begin{equation}\label{eq:extremal-root-form}
 \widetilde\mu_X\lrcorner\alpha_t =c(t)\overline{\alpha_t},\qquad |c(t)|=\rho(t)>0,
\end{equation}
with $c(t)$ constant on each component.

Fix a connected component $Z_{t,1}$ of $Z_t$ of degree $\ell$ over $X_t$. The $\nu$ cyclic group $\mathbb Z/\nu\mathbb Z$ acts on $Z_t$ and we let $H_1\leq \mathbb Z/\nu\mathbb Z$ to be the stabilizer group of $Z_{t,1}$ of order $\ell$. A generator of $H_1$ acts on $\alpha_t$ by a primitive $\ell$-th root $\zeta$ and preserves $\widetilde\mu_X$.
Equation \eqref{eq:extremal-root-form} gives $\zeta=\bar\zeta$, so $\ell\mid \gcdop(\nu,2)$.

Let $a:=\gcdop(\nu,2)$. On $Z_{t,1}$, the $a$-th power of the tautological section $\xi_t$ of $\pi_t^*K_{X_t}$ is invariant and descends to a holomorphic $a$-differential on $X_t$ via $\pi_t$.
Its $(\nu/a)$-th power is $s_t$.
Lemma~\ref{lem:saturation} therefore gives
\begin{align*}
 L_\nu=R^{\otimes (\nu/a)},\qquad R\subset E_a,\qquad \deg R=ad/2.
\end{align*}
By Theorem~\ref{thm:strict}, $R=L_a$, which proves assertions (1) and (2).

By \eqref{eq:extremal-equality-chain}, the equality $\rho_T=\rho$ on $U$ implies that the image of $U$ in $\mathcal M_g$ is a Teichm\"uller curve \cite[Theorem~5 and Corollary~1]{EKK94}.
\end{proof}

\subsection{A refined Arakelov inequality}
Let $\rho(t)^2|dt|^2$ be the Poincar\'e metric on $U=B\setminus D$ and let $\rho_T$ be the pullback Teichm\"uller metric defined in \eqref{eq:pullbackTeichmuller}. By the Schwarz lemma, $\rho_T\le\rho$. Hence, by the Gauss--Bonnet theorem,
\begin{equation}\label{eq:energy-area}
 \mathcal E_T(f):=\frac i2\int_U\rho_T^2(t)dt\wedge d\bar t\le\frac i2\int_U\rho^2dt\wedge d\bar t
 =\frac\pi2d.
\end{equation}
We prove the following refined Arakelov inequality for family of curves of genus $g\geq2$.

\begin{theorem}\label{thm:energy}
For every $\nu\ge1$ and every nonzero subsheaf $W\subseteq E_\nu$,
\begin{equation}\label{eq:energy-bound}
 \mu(W)\le\frac\nu\pi\mathcal E_T(f).
\end{equation}
If the classifying image of $U$ in $\mathcal M_g$ is not a Teichm\"uller curve, then we have
\[
 \qquad \mu_{\max}(E_\nu)\le\nu\left(\frac d2-\delta_f\right),
\]
where $\delta_f:=\frac d2-\frac{\mathcal E_T(f)}\pi>0$.
\end{theorem}
\begin{proof}
Replace $W$ by its saturation, of rank $r$, and let
$p:\PP(W)\to B$ parameterize lines. Write
$\xi=c_1(\OO_{\PP(W)}(1))$ and use $\gamma$ from
Proposition~\ref{prop:metric}. Choose smooth compactly supported forms
$0\le\gamma_j\le\gamma$ on $U$ whose masses approach
\[
 \int_U\gamma=d-\frac2\pi\mathcal E_T(f).
\]
They exist because $\gamma$ is continuous, nonnegative and integrable.
By Proposition~\ref{prop:metric}, the classes
\[
 A_j=2\xi+\nu p^*\bigl(c_1(\Om)-[\gamma_j]\bigr)
\]
are nef. Taking top intersections gives
\[
 0\le A_j^r
 =2^{r-1}\left(r\nu\left(d-\int_U\gamma_j\right)-2\deg W\right).
\]
Letting $j\to\infty$ proves \eqref{eq:energy-bound}.

If $\rho_T=\rho$ at one point, then the image of $U$ in $\mathcal M_g$ is a Teichm\"uller curve \cite[Theorem~5 and Corollary~1]{EKK94}. Otherwise $\rho_T<\rho$ everywhere. It follows that \eqref{eq:energy-area} is strict, and hence $\delta_f>0$.
\end{proof}

\section{Canonical equality in higher fibre dimension}\label{sec:geodesic}
We prove Theorem~\ref{thm:geodesic} by establishing a sharp Schwarz-type inequality between the Hodge metric/distance and the Poincar\'e metric/distance of any weight-$n$ $\mathbb C$PVHS on a unit disk $\Delta$ in Lemma~\ref{lem:block} and Lemma~\ref{lem:flatnorm}. Those bounds are sharp on the \(n\)-th symmetric power of the uniformizing variation of $\Delta$. For $L_1\subset E_1$ attains the Arakelov equality, we again use the Viehweg--Zuo iteration argument to construct maximal Higgs subbundle whose Hodge distance is exactly the Poincar\'e distance of $U$. This together with the previous distance inequality proves Theorem~\ref{thm:geodesic}.

\subsection{A sharp Hodge-norm estimate on a disk}\label{sec:disk}

Let $\Delta=\{z\in\mathbb C:|z|<1\}$ be the unit disk, equipped with the Poincar\'e metric of curvature $-4$ given by
\[
 ds_\Delta=\rho(z)|dz|,\qquad \rho(z)=\frac{1}{1-|z|^2}.
\]
Let $(E,\theta)$ be the graded Higgs bundle of an effective weight-$n$ $\mathbb C$-PVHS on $\Delta$, and write
\[
 E=\bigoplus_{k=0}^n E^{k,n-k},\qquad \theta=\bigoplus_{k=1}^n\theta^{k,n-k}\quad\text{ with }\quad \theta^{k,n-k}:E^{k,n-k}\longrightarrow E^{k-1,n-k+1}\otimes\Omega_\Delta^1.
\]
All norms below are taken with respect to the Hodge metric $h$.

\begin{lemma}\label{lem:block}
For every $z\in\Delta$ and $1\le k\le n$, we have
\begin{equation}\label{eq:block}
 \bigl\|\theta^{k,n-k}(\partial_z)\bigr\|_{\op}^2\le k(n-k+1)\rho(z)^2.
\end{equation}
\end{lemma}

\begin{proof}
Set $A_k(z):=\theta^{k,n-k}(\partial_z)$ and $f_k:=\|A_k\|_{\op}^2$.
Let $\Theta_k$ denote the Chern curvature of the Hodge metric $h_k$ on $E^{k,n-k}$ under the frame $dz\wedge d\bar z$.
The Griffiths' curvature formula; see \cite[Lemma~(7.18)]{Schmid}, gives
\begin{equation}\label{eq:harmonic-curvature}
 \Theta_k=A_k^*A_k-A_{k+1}A_{k+1}^*.
\end{equation}
At a point $z$ where $f_k>0$, choose unit norm vectors $u\in E_z^{k,n-k}$ and $v\in E_z^{k-1,n-k+1}$ satisfying
\[
 A_ku=\sqrt{f_k}\,v,\qquad A_k^*v=\sqrt{f_k}\,u.
\]
Extend $u\in E_z^{k,n-k}$ and the Hodge metric adjoint $v^*\in (E_z^{k-1,n-k+1})^\vee$ as local holomorphic sections such that $D_h u(z)=0$ and $D_h v^*(z)=0$ for the Chern connection $D_h$ of the Hodge metric $h=\oplus h_k$.
The function
\[
 F:=\frac{|v^*(A_ku)|^2}{\|u\|^2\|v^*\|^2}
\]
is smooth and positive near $z\in\Delta$, satisfies $F\le f_k$, and $F(z)=f_k(z)$.
Then \eqref{eq:harmonic-curvature} gives at the chosen point $z$
\begin{equation}\label{eq:support}
 \partial_z\partial_{\bar z}\log F=\langle\Theta_ku,u\rangle-\langle\Theta_{k-1}v,v\rangle=2f_k-\|A_{k+1}^*u\|^2-\|A_{k-1}v\|^2\ge 2f_k-f_{k-1}-f_{k+1}.
\end{equation}

Fix $0<R<1$ and put
\[
 \rho_R(z):=\frac{R}{R^2-|z|^2},\qquad c_k:=k(n-k+1).
\]
The functions $f_k$ are bounded on the closed disk of radius $R$, whereas $\rho_R$ diverges at its boundary.
Hence the maximum
\[
 M_R=\max_{\substack{|z|<R\\1\le k\le n}}\frac{f_k(z)}{c_k\rho_R(z)^2}
\]
is attained in the interior whenever some $f_k$ is nonzero.
Suppose $M_R>1$ and choose a pair $(z,k)$ such that $f_k(z)/c_k\rho_R^2(z)$ attains this maximum $M_R$.
Then the function $\log(F/(c_k\rho_R^2))$ has a local maximum there.
Since
\[
 \partial_z\partial_{\bar z}\log\rho_R^2=2\rho_R^2,\qquad 2c_k-c_{k-1}-c_{k+1}=2,
\]
equation \eqref{eq:support} gives
\[
 0\ge \partial_z\partial_{\bar z}\log \bigl(F(z)/c_k\rho_R(z)^2\bigr)\ge 2f_k-f_{k-1}-f_{k+1}-2\rho_R^2\ge 2(M_R-1)\rho_R^2>0,
\]
a contradiction.
Therefore $M_R\le 1$, and letting $R$ increase to $1$ proves \eqref{eq:block}.
\end{proof}

\begin{lemma}\label{lem:flatnorm}
Let $u$ be a nonzero flat section of the weight-$n$ $\mathbb C$-PVHS on $\Delta$.
Then
\begin{equation}\label{eq:flatnorm}
 \left|\log\frac{\|u(w)\|}{\|u(z)\|}\right|\le n\,d_\Delta(z,w),\qquad\text{for any}\quad z,w\in\Delta,
\end{equation}
where $\|u\|$ is the Hodge norm and $d_\Delta$ is the Poincar\'e distance.
\end{lemma}

\begin{proof}
For a real tangent vector $\xi$ of $\Delta$ at $z$, let $\Phi(\xi):=\theta(\xi^{1,0})+\theta(\xi^{1,0})^*$.
This endomorphism is self-adjoint.
Write $e=\bigoplus_{k=0}^n e^{k,n-k}\in E_z$ and put $c_k=k(n-k+1)$ for $1\le k\le n$.
By \eqref{eq:block},
\[
 |\langle\Phi(\xi)e,e\rangle|\le 2\rho(z)|dz(\xi)|\sum_{k=1}^n\sqrt{c_k}\,\|e^{k-1,n-k+1}\|\cdot\|e^{k,n-k}\|.
\]
Let $J$ be the $(n+1)\times(n+1)$ real symmetric matrix with $J_{k-1,k}=J_{k,k-1}=\sqrt{c_k}$ for $1\leq k\leq n$ and $J_{i,j}=0$ otherwise.
Let $x:=(\sqrt{\binom n0},\sqrt{\binom n1},\cdots,\sqrt{\binom nn})\in\mathbb R^{n+1}$. Then $x$ is an eigenvector of $J$ with eigenvalue $n$.
By the Perron--Frobenius theorem, $n$ is the largest eigenvalue of $J$.
Since \[
2\sum_{k=1}^n\sqrt{c_k}\,\|e^{k-1,n-k+1}\|\cdot\|e^{k,n-k}\|=(\|e^{0,n}\|,\cdots,\|e^{n,0}\|)\,J\,(\|e^{0,n}\|,\cdots,\|e^{n,0}\|)^{\mathsf T},\]
it follows that
\begin{equation}\label{eq:realhiggs}
 \|\Phi(\xi)\|_{\op}\le n\rho(z)|dz(\xi)|.
\end{equation}

The smooth flat connection is $\nabla=D_h+\theta+\theta^*$, where $D_h$ is the Chern connection of the Hodge metric.
Since $u$ is flat, $D_{h,\xi}u=-\Phi(\xi)u$, and hence
\[
 |d\log\|u\|(\xi)|=\left|\frac{\langle\Phi(\xi)u,u\rangle}{\|u\|^2}\right|\leq n\rho(z)|dz(\xi)|,
\]
by \eqref{eq:realhiggs}.
Integrating along the real Poincar\'e geodesic from $z$ to $w$ proves \eqref{eq:flatnorm}.
\end{proof}

\subsection{Canonical Arakelov equality and Kobayashi geodesicity}

\begin{proof}[Proof of Theorem~\ref{thm:geodesic}]
Let $D$ be the divisor of the singular fiber and let $f^\circ:=f|_U:X^\circ\to U:=B\setminus D$ be the family of smooth varieties. By \cite[Proposition~1.1(c)]{VZ06}, the maximal degree line bundle $L\subset E_1$ generates a sub $\mathbb C$PVHS $\mathbb W$ of the middle-cohomology PVHS of $f^\circ$. And the logarithmic Higgs bundle of $\mathbb W$ is
\begin{align*}
W&=\bigoplus_{k=0}^n W^{k,n-k},\qquad W^{k,n-k}\simeq L\otimes\Om^{k-n},\\ \theta^{k,n-k}:&W^{k,n-k}\xrightarrow{\ \simeq\ }W^{k-1,n-k+1}\otimes\Om\quad(1\le k\le n).\end{align*}

By \cite[Proposition~3.4]{VZ03}, applied possibly after a finite \'etale cover of $U$, the pullback of $(W,\theta)$ to $\Delta$ is the $n$-th symmetric power of the rank-two uniformizing Higgs bundle of $\Delta$, up to a unitary twist.

Now we prove the classifying map $\widetilde\varphi:\Delta\to \widetilde M$ is a Kobayashi isometry. Fix $z,w\in\Delta$ and set $a:=d_\Delta(z,w)$. After taking an automorphism of $\Delta$ in $\mathrm{SU}(1,1)$, we may take $z=0$ and $w=t:=\tanh a$.
Let $(u_1,u_2)$ be a flat basis of the uniformization Higgs bundle of $\Delta$ such that its Hodge metric $h$ at $0\in\Delta$ with respect to this flat basis is $\operatorname{diag}\,(1,1)$, the element
\[
g_t=\frac{1}{\sqrt{1-t^2}}\begin{pmatrix}1&t\\t&1\end{pmatrix}\in\mathrm{SU}(1,1)
\]
acts on $\Delta$ by $\zeta\mapsto(\zeta+t)/(1+t\zeta)$ and sends $0$ to $t$.
Since $g_t$ preserves the polarization, we have $h_t(g_t v_1,g_t v_2)=h_0(v_1,v_2)$ for any vectors $v_1,v_2$ at $0$. Thus, in the basis $(u_1,u_2)$, the matrix of $h_t$ is
\[
(g_t^{-1})^*g_t^{-1}=\rho(t)\begin{pmatrix}1+t^2&-2t\\-2t&1+t^2\end{pmatrix}.
\]
The eigenvalues are $e^{2a}$ and $e^{-2a}$. By taking $n$-th power of the uniformization Higgs bundle, there is a nonzero flat vector $u$ of $\mathbb W$ on $\Delta$ such that
\begin{equation}\label{eq:endpoint}
 \left|\log\frac{\|u(\widetilde\varphi(w))\|}{\|u(\widetilde\varphi(z))\|}\right|=na=nd_\Delta(z,w).
\end{equation}

On $\widetilde{\mathcal M}$, $u$ extends uniquely to a global flat section of the pullback middle-cohomology VHS of $f^\circ$. We denote this section again by $u$.
The extended section still satisfies the identity \eqref{eq:endpoint}.
By definition, the Kobayashi pseudodistance is
\[
d_{\widetilde{\mathcal M}}\bigl(\widetilde\varphi(z),\widetilde\varphi(w)\bigr)=\inf\sum_{i=1}^N d_\Delta(a_i,b_i),
\]
where the infimum is over all finite chains of holomorphic maps $\psi_i:\Delta\to\widetilde{\mathcal M}$ and points $a_i,b_i\in\Delta$ satisfying $\psi_1(a_1)=\widetilde\varphi(z)$, $\psi_N(b_N)=\widetilde\varphi(w)$, and $\psi_i(b_i)=\psi_{i+1}(a_{i+1})$ for $1\le i<N$.

For each map $\psi_i$, $\psi_i^*u$ is a flat section of the pullback $\mathbb C$PVHS. Lemma~\ref{lem:flatnorm} therefore gives
\[
\left|\log\frac{\|u(\psi_i(b_i))\|}{\|u(\psi_i(a_i))\|}\right|\le n\,d_\Delta(a_i,b_i).
\]
Using \eqref{eq:endpoint} and the triangle inequality, we obtain
\[
n\,d_\Delta(z,w)=\left|\sum_{i=1}^N\log\frac{\|u(\psi_i(b_i))\|}{\|u(\psi_i(a_i))\|}\right|\le n\sum_{i=1}^N d_\Delta(a_i,b_i).
\]
Taking the infimum over all such chains gives
\begin{align*}
n\,d_\Delta(z,w)\le n\,d_{\widetilde{\mathcal M}}\bigl(\widetilde\varphi(z),\widetilde\varphi(w)\bigr).
\end{align*}
Since $\widetilde\varphi$ is holomorphic, the distance-decreasing property of the Kobayashi distance gives the reverse inequality.
Therefore $\widetilde\varphi$ is a Kobayashi isometric embedding.
\end{proof}

\section{The number of singular fibres and its consequences}
\label{sec:surfaces}
We return to non-isotrivial semistable curve fibrations of genus
$g\ge2$. The main assertion over $\PP^1$ uses only the slope bounds
and the curve equality classification already proved.

\subsection{The singular-fibre theorem and the six-fibre ring}
By Theorem~\ref{thm:strict} with $\nu=1$ and $W=E_1$, we have
Tan's strict Hodge inequality \eqref{eq:LiuTan}.
Beauville proved $s\ge4$, and in the four-fibre case
$p_g(X)=q(X)=0$ \cite{Beauville81}. Then $\deg E_1=g$ and $d=2$,
contradicting \eqref{eq:LiuTan}. Thus $s\ge5$, recovering Tan's theorem
\cite{Tan95}.
Write $P_\nu(X)=h^0(X,\omega_X^{\otimes\nu})$ and
$R(X,K_X)=\bigoplus_{\nu\ge0}H^0(X,\omega_X^{\otimes\nu})$.

\begin{theorem}\label{thm:six}
Suppose that $B=\PP^1$, and put $s=|D|$.
\begin{enumerate}
\item If $s\le5$, then $P_\nu(X)=0$ for every $\nu\ge1$, and $\kappa(X)=-\infty$.
\item If $s=6$, then $P_\nu(X)\le1$ for every $\nu\ge1$, and $\kappa(X)\le0$.
\item If $\kappa(X)\ge1$, then $s\ge7$.
\end{enumerate}
\end{theorem}
\begin{proof}
Since $\omega_{\PP^1}=\OO_{\PP^1}(-2)$, it follows from the projection formula that
\begin{equation}\label{eq:plurigenus}
 P_\nu(X)=\dim\Hom\bigl(\OO_{\PP^1}(2\nu),E_\nu\bigr).
\end{equation}
If $s\le5$, then $\nu d/2=\nu(s-2)/2<2\nu$. A nonzero map in \eqref{eq:plurigenus} would have a saturated line image of degree at least $2\nu$, contradicting the line bound. Thus $P_\nu(X)=0$ for every $\nu\ge1$.

If $s=6$, then $d=4$. Corollary~\ref{coro:Hom}, applied to
$M=\OO_{\PP^1}(2\nu)$, gives $P_\nu(X)\le1$ for every $\nu\ge1$.

If the Kodaira dimension is positive, some pluricanonical system has positive-dimensional image and hence at least two independent sections. This contradicts the preceding bounds when $s\le6$, and proves (3).
\end{proof}

Tan--Tu--Zamora proved that $\kappa(X)\ge0$ implies $s\ge6$ \cite{TTZ05}. Lu--Tan--Zuo proved $s\ge7$ for surfaces of general type and described possible surfaces with $\kappa(X)=1$, $s=6$ \cite[Theorems~1 and~3]{LTZ18}. In their Theorem~3 and the preceding paragraph, Lu--Tan--Zuo left the existence of the latter case undecided. By Theorem~\ref{thm:six}, this case cannot occur. In particular, the elliptic surfaces with multiple-fibre multiplicities $(2,3)$ and $(2,5)$ listed there cannot be minimal models of a semi-stable fibration with exactly six singular fibers over $\mathbb{P}^1$. The general-type bound is also recalled in \cite{CMLZ24,Zamora24}.

\begin{remark}
One can also use the splitting $E_\nu=\bigoplus_i\OO(a_i)$ on $\PP^1$. When $s=6$, every $a_i\le2\nu$, and at most one $a_i$ equals $2\nu$. Thus $E_\nu(-2\nu)$ has at most one trivial summand, and all other summands have negative degree.
\end{remark}

We next apply the parity theorem to the canonical ring.

\begin{theorem}\label{thm:ring}
If $B=\PP^1$ and $s=6$, then the graded canonical ring is one of the following three possibilities.
\[
 \CC,\qquad \CC[u]\quad(\deg u=1),\qquad
 \CC[v]\quad(\deg v=2).
\]
When $\kappa(X)=0$, the marked classifying map identifies the universal cover of $U$ with a quadratic Teichm\"uller disk whose effective stabilizer is a lattice. For such a surface, the maximal line comes from a holomorphic one-form if $p_g(X)>0$, and from a nonsquare holomorphic quadratic differential on a general fibre if $p_g(X)=0$.
\end{theorem}
\begin{proof}
Here $d=s-2=4$. By \eqref{eq:plurigenus}, $P_\nu(X)=h^0(\PP^1,E_\nu(-2\nu))$.
Every summand in the splitting of $E_\nu$ has degree at most $2\nu$, so for every $\nu\ge1$,
\begin{equation}\label{eq:detect}
	P_\nu(X)>0\quad\Longleftrightarrow\quad
	\mu_{\max}(E_\nu)=2\nu.
\end{equation}
By Theorem~\ref{thm:six}, $P_\nu(X)\le1$ for all $\nu\ge1$. If the equality in the Arakelov bound occurs at no exponent, \eqref{eq:detect} therefore gives $R(X,K_X)=\CC$.

Suppose the equality occurs at some odd exponent $\nu$, and let $L_\nu\subset E_\nu$ be the maximal line.
By Theorem~\ref{thm:characterization}(1), there is a saturated line $A\subset E_1$ of degree $d/2=2$ such that
$L_\nu=A^{\otimes\nu}$ under multiplication of differentials.
Thus $A\simeq\OO(2)$ gives a nonzero $u\in H^0(X,K_X)$. For every $m\ge1$, $u^m$ is nonzero and spans $H^0(X,mK_X)$. Hence $R(X,K_X)=\CC[u]$.

Now suppose the equality occurs at $\nu=2k$ and at no odd exponent. By Theorem~\ref{thm:characterization}(2),
the maximal line satisfies $L_{2k}=Q^{\otimes k}$ for a saturated line $Q\subset E_2$ of degree $d=4$.
Thus $Q\simeq\OO(4)$ gives a nonzero $v\in H^0(X,2K_X)$. Its powers span the even graded pieces, and the odd pieces vanish by \eqref{eq:detect}. Hence $R(X,K_X)=\CC[v]$.
If a generator of $Q$ on a general fibre $F$ were the square of a holomorphic one-form on $F$, Lemma~\ref{lem:saturation} would give a saturated line $A\subset E_1$ with $Q=A^{\otimes2}$. Then $\deg A=2=d/2$, so equality would occur at exponent one, contrary to the hypothesis.

Finally, if $\kappa(X)=0$, some $P_\nu(X)$ is positive. By \eqref{eq:detect} and the proof of Theorem~\ref{thm:characterization}, the marked classifying map identifies the universal cover of $U$ with a Teichm\"uller disk, whose effective stabilizer is a lattice. The two cases above distinguish $p_g(X)=1$ from $p_g(X)=0$.
\end{proof}

\begin{corollary}\label{cor:index}
Let $B=\PP^1$, $s=6$ and $\kappa(X)=0$. If $S$ is the minimal model of $X$ and $r$ is the order of its canonical bundle, then $r\mid2$. In particular, the canonical-index $3$, $4$ and $6$ cases cannot occur.
\end{corollary}
\begin{proof}
By Theorem~\ref{thm:ring}, $P_2(X)=1$. It follows from the birational invariance of plurigenera that $H^0(S,2K_S)\ne0$. A torsion line bundle on a projective surface with a nonzero section is trivial. Indeed, its effective zero divisor is numerically trivial and hence has zero intersection with an ample divisor. Therefore $2K_S\simeq\OO_S$.
\end{proof}

\begin{remark}
Lu--Tan--Zuo obtained the Teichm\"uller description for surfaces of Kodaira dimension zero after a suitable finite \emph{\'etale} cover of the total surface \cite[Theorem~2]{LTZ18}. Here the quadratic formulation applies to the original family. Tu constructed six-fibre examples on surfaces of Kodaira dimension zero \cite{Tu07}; this does not assert existence for every possible canonical index. By Corollary~\ref{cor:index}, we obtain a restriction on the possible types.
\end{remark}

Combining the six-fibre result with an effective quadratic differential, we now exclude maximal destabilization on surfaces of Kodaira dimension one.

\begin{lemma}\label{lem:horizontal-big}
Let $Z$ be an effective divisor on a projective fibred surface $f:X\to B$, and suppose that $Z\cdot F>0$ for a general fibre $F$. If $L$ is a line bundle of positive degree on $B$, then $Z+f^*L$ is big.
\end{lemma}
\begin{proof}
Fix an ample divisor $H$ on $X$. By Riemann--Roch on the generic fibre, we can choose $m_0>0$ such that the restriction of $m_0Z-H$ has a nonzero section. Twisting its nonzero direct image by a sufficiently positive line bundle $A$ on $B$, we obtain an effective divisor $m_0Z+f^*A-H$. Choose a positive rational $\epsilon$ such that $\epsilon m_0<1$ and $\deg(L-\epsilon A)>0$. Then
\[
 Z+f^*L=\epsilon H+
 \epsilon(m_0Z+f^*A-H)+(1-\epsilon m_0)Z+f^*(L-\epsilon A).
\]
After taking a sufficiently divisible multiple, the terms after $\epsilon H$ are effective as rational divisor classes. Hence the divisor is big.
\end{proof}

\begin{theorem}\label{thm:Kodaira}
Suppose that a pluricanonical direct image attains its Arakelov bound at an exponent $\nu\ge2$. Write $b=g(B)$ and $s=|D|$. Then
\[
 \kappa(X)=
 \begin{cases}
 -\infty,&6b-6+s<0,\\
 0,&6b-6+s=0,\\
 2,&6b-6+s>0.
 \end{cases}
\]
Consequently, if $\kappa(X)=1$, then $\mu_{\max}(E_\nu)<\nu d/2$ for every $\nu\ge2$.
\end{theorem}
\begin{proof}
By the parity theorem, there is a maximal quadratic line $Q\subset E_2$, obtained from $A^{\otimes2}$ in the odd case. Since $Q\simeq\Om$, evaluation gives an effective divisor
\begin{equation}\label{eq:critical-divisor}
 Z\in|2K_f-f^*(K_B+D)|=|2K_X-f^*(3K_B+D)|,
 \qquad Z\cdot F=4g-4>0.
\end{equation}
If $\deg(3K_B+D)=6b-6+s>0$, then $2K_X$ is big by Lemma~\ref{lem:horizontal-big}.

Since $d=2b-2+s>0$, the remaining cases have $b=0$. A negative value means $s\le5$. By Theorem~\ref{thm:six}, we have $\kappa(X)=-\infty$ in this case. A zero value means $s=6$, and a maximal line $L_\nu\simeq\OO(2\nu)$ gives a nonzero pluricanonical section. On the other hand, $\kappa(X)\le0$ by Theorem~\ref{thm:six}. Hence $\kappa(X)=0$. This proves the three cases and excludes equality when $\kappa(X)=1$. Equality at $\nu=1$ is also impossible, since equality at $\nu=2$ would then follow from Lemma~\ref{lem:saturation}.
\end{proof}

\begin{proof}[Proof of Theorem~\ref{thm:intro-surfaces}]
The assertions over $\PP^1$ follow from Theorems~\ref{thm:six}
and~\ref{thm:ring}. The final assertion follows from
Theorem~\ref{thm:Kodaira}, including degree one by
Lemma~\ref{lem:saturation}.
\end{proof}

\subsection{Another proof of Tan's strict canonical class inequality}

Assume that $f$ is relatively minimal. Put $K_f=K_X-f^*K_B$ and $\chi_f=\deg E_1$. By the Riemann--Roch formula and relative duality, we have
\begin{equation}\label{eq:GRR}
 \deg E_\nu=\frac{\nu(\nu-1)}2K_f^2+\chi_f,
 \qquad \rk E_\nu=(2\nu-1)(g-1)\quad(\nu\ge2).
\end{equation}
%Indeed, $\omega_{X_b}$ has nonnegative degree on every component and positive degree on at least one. Therefore, a section of $\omega_{X_b}^{1-\nu}$ vanishes on every positive-degree component. By connectedness it also vanishes on the degree-zero components. Hence $R^1f_*\omega_{X/B}^{\otimes\nu}=0$. For a nonminimal nodal model, the right side of the degree formula has the additional term $\operatorname{length}(R^1f_*\omega_{X/B}^{\otimes\nu})\ge0$.

Apply \eqref{eq:energy-bound} to the full bundle $E_\nu$. Substituting \eqref{eq:GRR}, dividing by $\nu^2$, and letting $\nu\to\infty$, we obtain
\begin{equation}\label{eq:energy-K2}
 K_f^2\le\frac{4(g-1)}\pi\mathcal E_T(f).
\end{equation}
We distinguish two cases according to whether the classifying image of $U$ in $\mathcal M_g$ is a Teichm\"uller curve.

\medskip
\noindent\emph{The non-Teichm\"uller case.}
Thus $\mathcal E_T(f)<\pi d/2$, and \eqref{eq:energy-K2} gives $K_f^2<2(g-1)d$.

\medskip
\noindent\emph{The Teichm\"uller case.}
Let $m_1,\ldots,m_\ell$ be the positive zero orders of its generating holomorphic quadratic differential, so $\sum_i m_i=4g-4$.
The intersection formula of Chen--M\"oller gives
\begin{equation}\label{eq:CM}
 \frac{K_f^2}{d}=\frac14\sum_i\frac{m_i(m_i+4)}{m_i+2}
 =g-1+\frac12\sum_i\frac{m_i}{m_i+2}
\end{equation}
by \cite[Proposition~4.2]{CM14} for nonsquare differentials and \cite[Proposition~4.8]{CM12} for squares.
Hence
\begin{equation}\label{eq:five-thirds}
 K_f^2\le\frac53(g-1)d.
\end{equation}
Equality holds precisely when all zeros are simple.
If the quadratic differential is a square, write $m_i=2a_i$, so that $\sum_i a_i=2g-2$. In this case the bound improves to
\[
 K_f^2\le\frac32(g-1)d.
\]
Equality holds precisely when the one-form has only simple zeros. The latter bound also appears in \cite[Proposition~1]{LTZ18}.

The boundary conditions in \eqref{eq:CM} follow from \cite[Corollaries~4.5--4.6]{CM14}. After a finite cover unramified over $U$, the zero sections remain distinct and avoid the nodes, and the differential is nonzero on every component of the stable limit. For the Abelian case, see \cite[Proposition~3.1 and Corollary~3.2]{CM12}.
On the minimal resolution, write $Z=\sum_i m_iS_i=2K_f-f^*(K_B+D)$. Then
\[
 S_i^2=-\frac d{m_i+2},\qquad Z^2=4K_f^2-8(g-1)d.
\]
Further vertical blowups require the corresponding exceptional-divisor corrections.

Together, \eqref{eq:energy-K2} and \eqref{eq:five-thirds} give another proof of Tan's strict canonical class inequality
\[
 K_f^2<(2g-2)d.
\]

\subsection{Splitting over the projective line}\label{sec:splitting}

Assume here that $f$ is relatively minimal, let $B=\PP^1$ and $\nu\ge2$. Put
\[
 N_\nu=(2\nu-1)(g-1),\qquad c_\nu=\frac{\nu(s-2)}2,
 \qquad E_\nu=\bigoplus_{i=1}^{N_\nu}\OO(a_i),
\]
with $a_1\ge\cdots\ge a_{N_\nu}$.

\begin{proposition}\label{prop:splitting}
If $c_\nu$ is an integer, then
\[
 a_1\le c_\nu,\quad a_2\le c_\nu-1,\qquad
 \deg E_\nu\le N_\nu(c_\nu-1)+1.
\]
If $c_\nu$ is a half-integer, then $\deg E_\nu\le N_\nu\lfloor c_\nu\rfloor$. If no maximal line occurs and $c_\nu$ is integral, then
$\deg E_\nu\le N_\nu(c_\nu-1)$.
\end{proposition}
\begin{proof}
By the line bound, every summand satisfies $a_i\le c_\nu$. Two summands of degree $c_\nu$ would give a rank-two equality subbundle, contrary to the strict rank-two inequality. Summing these integral bounds, we obtain the degree estimates.
\end{proof}

For a relatively minimal fibration, combining \eqref{eq:GRR} with Proposition~\ref{prop:splitting} at $\nu=2$, we obtain
\begin{equation}\label{eq:geography}
 K_f^2+\chi_f\le3(g-1)(s-3)+1.
\end{equation}
The final $1$ can be removed if there is no maximal line. If equality holds in \eqref{eq:geography}, then
\[
 E_2\simeq\OO(s-2)\oplus\OO(s-3)^{\oplus(3g-4)}
\]
so the marked classifying map uniformizes a Teichm\"uller disk.

\begin{remark}
For a nonminimal nodal model, the degree formula has an additional nonnegative term:
\[
 \deg E_\nu=\frac{\nu(\nu-1)}2K_f^2+\chi_f+
 \operatorname{length}\bigl(R^1f_*\omega_{X/B}^{\otimes\nu}\bigr).
\]
Therefore \eqref{eq:geography} remains valid, but the splitting in the equality case is asserted here only for a relatively minimal model. We do not use this formula in the proof of Theorem~\ref{thm:six}. By taking the limit in the unrefined estimate $\mu(E_\nu)<\nu d/2$, we obtain only a non-strict canonical class bound. By Theorem~\ref{thm:energy}, we obtain the additional metric gap used above.
\end{remark}

\raggedbottom

\end{document}